\documentclass[12pt, one side,emlines]{amsart}
\usepackage[a4paper,margin=1in]{geometry}
\usepackage{fancyhdr}
\usepackage{enumerate}
\usepackage{amssymb,latexsym,xy,eucal,mathrsfs}
\textwidth=18cm \textheight=27cm \theoremstyle{plain}
\usepackage{tikz}
\usetikzlibrary{positioning}
\theoremstyle{definition}

\newtheorem{theorem}{Theorem}[section]
\newtheorem{thm}[theorem]{Theorem}
\newtheorem{lem}[theorem]{Lemma}

\newtheorem{prop}[theorem]{Proposition}

\newtheorem{defn}{Definition}[section]

\begin{document}
	

	%
	\title{Cographic Splitting of Graphic Matroids with Respect to a Set with  Three Elements}\maketitle
	
	\markboth{ Shital Dilip Solanki, Ganesh Mundhe and S. B. Dhotre}{ Cographic Splitting of Graphic Matroid with Respect to Three Elements }\begin{center}\begin{large} Shital Dilip Solanki$^1$, Ganesh Mundhe$^2$ and S. B. Dhotre$^3$ \end{large}\\\begin{small}\vskip.1in\emph{
				1. Ajeenkya DY Patil University, Pune-411047, Maharashtra,
				India\\ 
				2. Army Institute of Technology, Pune-411015, Maharashtra,
				India\\
				3. Department of Mathematics,
				Savitribai Phule Pune University,\\ Pune - 411007, Maharashtra,
				India}\\
			E-mail: \texttt{1. shital.solanki@adypu.edu.in, 2. gmundhe@aitpune.edu.in, 3. dsantosh2@yahoo.co.in. }\end{small}\end{center}\vskip.2in
	\begin{abstract} In general, the splitting operation on binary matroids does not preserve the graphicness and cographicness properties of binary matroids. In this paper, we obtain a characterization of the class of graphic matroids whose splitting with respect to a set of three elements gives cographic matroids. We also provide an alternate proof to such characterization when the set contains two elements which was provided by Borse et al. \cite{ymb1}.	\end{abstract}\vskip.2in
	
	\noindent\begin{Small}\textbf{Mathematics Subject Classification (2010)}:
		05B35, 05C50, 05C83     \end{Small}\vskip.2in
	\vskip.25in

	\baselineskip 19truept 
	\section{Introduction}
	\noindent We refer to Oxley \cite{ox} for undefined terminology and notations. There are many operations in graph theory that are of great importance. The Splitting operation is one of them, which is introduced by Fleischner \cite{fl}. He has characterized Eulerian graphs using this operation. The Splitting operation is given as follows.
	Let $v$ be a vertex of a connected graph G such that $d(v)\geq 3$ . If $x,y$ are two edges incident at $v$ such that $x=v v_{1}$ and $y=v v_{2}$. 
	Let $G_{x, y}$ be obtained from $G$ by removing edges $x,y$ and by adding a vertex $v_{x,y}$ adjacent to $v_{1}$ and $v_{2}$. The splitting operation on the graph $G$ with respect to $x$ and $y$ is defined as the transition from $G$ to $G_{x, y}$.

	The Splitting operation is further extended to binary matroids by Raghunathan et al. \cite{ttr} and defined the splitting operation with respect to a pair of elements of a binary matroid. Later, it is generalized by Shikare at al. \cite{mms1} by defining the splitting of a binary matroid with respect to any set $T$ as follows. 
	\begin{defn}\cite{mms1}\label{spdef}
		Let $M$ be a binary matroid and let $A$ be the standard matrix representation of $M$ over the field $GF(2)$ and let $T \subseteq E(M)$. Obtain a matrix $A_{T}$ by adjoining one row at the bottom of the matrix $A$ with entries 1 in the columns corresponding to the elements of the set $T$ and zero otherwise. Then the vector matroid, denoted by $M_T$, is called  the splitting matroid of $M$ with respect to $T$, and the splitting operation with respect to $T$ is the transition from $M$ to $M_{T}$.
	\end{defn}
	Some properties of a binary matroid are not preserved under the splitting operation. Connectedness, graphicness and cographicness are a few properties among the properties of binary matroids which are not preserves under the splitting operation. 
	Shikare and Waphare \cite{mms} obtained characterization for a graphic matroid $M$ whose splitting matroid $M_T$ is also graphic for any $T$ with two elements. Mundhe \cite{gm} obtained a similar characterization for a set $T$ with three elements. Borse et al. \cite{ymb1} characterized the class of graphic matroids, which yields a cographic matroid under the splitting operation with respect to a pair of elements. In this paper, we characterize a graphic matroid that gives a cographic matroid under the splitting operation with respect to any set $T$ with $|T|=3$.

	Let $\mathcal{M}_k$ be the class of graphic matroid whose splitting with respect to any set of $k$-elements is cographic. Let $N$ be an extension of $M(K_5)$ or $M(K_{3,3})$ by an element $a$. Then the splitting of $N$ with respect to $a$ is not cographic. Therefore the single element extensions of $M(K_5)$ and $M(K_{3,3})$  are the forbidden minors the class $\mathcal{M}_1$ and only these are the forbidden minors as per the definition of the splitting operation. 
	Let $N$ be an extension of $M(K_5)$  by a set $T$ with $|T|=2$  such that $T$ be a 2-cocircuit of $N$. Then $N\backslash T= M(K_5)$ and, by Definition \ref{spdef}, $N_T\backslash T=N\backslash T=M(K_5).$ Thus $N$ is a forbidden minor for the class $\mathcal{M}_2.$ Similarly,  it is easy to check that $\mathcal{M}_2$ contains  a forbidden minor $N$ which arise as follows. (i) $N$ is an extension of $M(K_5)$ or  $M(K_{3,3})$ by two elements $x,y$ such that ${x,y}$ is a 2-cocircuit of $N$ or (ii) $N$ is a coextension of $M(K_5)$ or  $M(K_{3,3})$ by an element $x$ such that ${x,y}$ is a 2-cocircuit of $N$ for some $y \in E(N)$.

If $N$ is a forbidden minor for the class $\mathcal{M}_k$ which is isomorphic to $M(K_ 5) $,  $M(K_{3,3})$ or a matroid which arises by extensions and coextensions   of $M(K_ 5)$ or $M(K_{3,3}),$ then we call $N$ is a trivial forbidden minor for the class $\mathcal{M}_k$.  Let $\mathcal{T}_k$ be a set containing all trivial forbidden minors for the class $\mathcal{M}_k$. 
 We focus on the forbidden minors of the class $\mathcal{M}_k$ other than the members of $\mathcal{T}_k$.

   In the following result,  Borse et al. \cite{ymb1} have given forbidden minors for the class $\mathcal{M}_2.$
\begin{thm}\cite{ymb1}\label{gct}
	Let $M$ be a graphic matroid without containing a minor isomorphic to any member of $\mathcal{T}_2$. Then $M \in \mathcal{M}_2$ if and only if $M$ does not contain a minor isomorphic to any one of the circuit matroids $M(G_1)$, $M(G_2)$ and $M(G_3)$, where $G_1$, $G_2$ and $G_3$ are the graphs as shown in Figure \ref{1}.
\end{thm}
\begin{figure}[h!]
	\centering
\unitlength 1mm 
\linethickness{0.4pt}
\ifx\plotpoint\undefined\newsavebox{\plotpoint}\fi 
\begin{picture}(89.119,31.693)(0,0)
	\put(3.569,8.056){\circle*{1.682}}
	\put(21.438,8.024){\circle*{1.682}}
	\put(3.893,22.31){\circle*{1.682}}
	\put(21.762,22.278){\circle*{1.682}}
	\put(3.751,22.316){\line(1,0){18}}
	\put(3.677,22.312){\line(0,-1){14.125}}
	\put(21.552,7.937){\line(0,1){14.625}}
	\qbezier(3.569,22.561)(12.924,29.078)(21.859,22.141)
	\put(3.779,8.056){\line(1,0){18.079}}
	\qbezier(3.359,8.056)(13.24,.067)(21.438,8.056)
	\multiput(3.674,22.351)(.0418961088,-.0337152118){424}{\line(1,0){.0418961088}}
	\multiput(3.359,8.056)(.0433835446,.0337152118){424}{\line(1,0){.0433835446}}
	\put(28.16,7.824){\circle*{1.682}}
	\put(46.029,7.791){\circle*{1.682}}
	\put(28.484,22.078){\circle*{1.682}}
	\put(46.353,22.046){\circle*{1.682}}
	\put(28.342,22.084){\line(1,0){18}}
	\put(28.268,22.08){\line(0,-1){14.125}}
	\put(46.143,7.705){\line(0,1){14.625}}
	\put(28.37,7.824){\line(1,0){18.079}}
	\qbezier(27.95,7.824)(37.83,-.165)(46.029,7.824)
	\put(37.823,30.852){\circle*{1.682}}
	\multiput(37.874,31.052)(.0337022512,-.0361095549){247}{\line(0,-1){.0361095549}}
	\multiput(28.248,22.189)(.0364583607,.0336174495){264}{\line(1,0){.0364583607}}
	\multiput(27.998,7.939)(.033670059,.0778620114){297}{\line(0,1){.0778620114}}
	\put(61.028,8.421){\circle*{1.682}}
	\put(78.897,8.388){\circle*{1.682}}
	\put(61.352,22.675){\circle*{1.682}}
	\put(79.221,22.643){\circle*{1.682}}
	\put(61.21,22.681){\line(1,0){18}}
	\put(79.011,8.302){\line(0,1){14.625}}
	\put(61.238,8.421){\line(1,0){18.079}}
	\multiput(61.133,22.716)(.0418962264,-.0337146226){424}{\line(1,0){.0418962264}}
	\multiput(60.818,8.421)(.0433820755,.0337146226){424}{\line(1,0){.0433820755}}
	\put(53.028,16.398){\circle*{1.682}}
	\put(88.278,16.523){\circle*{1.682}}
	\multiput(52.687,16.432)(8.93751,.03125){4}{\line(1,0){8.93751}}
	\multiput(79.062,8.432)(.038381772,.033713718){241}{\line(1,0){.038381772}}
	\multiput(88.312,16.557)(-.050403264,.033602176){186}{\line(-1,0){.050403264}}
	\multiput(52.937,16.682)(.047857179,.033571454){175}{\line(1,0){.047857179}}
	\multiput(52.812,16.307)(.035326113,-.033695677){230}{\line(1,0){.035326113}}
	\put(13.258,.25){\makebox(0,0)[cc]{$G_1$}}
	\put(37.3,.25){\makebox(0,0)[cc]{$G_2$}}
	\put(68.766,.25){\makebox(0,0)[cc]{$G_3$}}
\end{picture}

	\caption{Forbidden Minors for the Class $\mathcal{M}_2$.}
	\label{1}
	
\end{figure}

Denote, by $\tilde{M}$, the single element extension of $M$. 
Later on,  Mundhe \cite{gm} obtained forbidden minor for a class of graphic matroid which yields graphic matroid under the splitting operation with respect to a set of three elements as follows.
\begin{thm}\cite{gm} 
	Let $M$ be a graphic matroid. Then $M_T$ is graphic for any  $T\subseteq E(M)$ with $|T|=3$ if and only if $M$ does not contain a minor isomorphic to any of the circuit matroids $\tilde{M}(Q_1)$, $\tilde{M}(Q_2)$,  $M(Q_3),$ $M(Q_4),$ $M(Q_5)$ and $M(Q_6),$ where $Q_i$ is the graph as shown in 
	Figure \ref{gtg}, for $i=1,2,\dots,6.$
\end{thm}
\begin{figure}[h!]
\centering
\unitlength 1mm 
\linethickness{0.4pt}
\ifx\plotpoint\undefined\newsavebox{\plotpoint}\fi 
\begin{picture}(162.141,34.875)(0,0)
	\put(76.591,10.721){\circle*{1.682}}
	\put(94.46,10.689){\circle*{1.682}}
	\put(76.915,24.975){\circle*{1.682}}
	\put(94.784,24.943){\circle*{1.682}}
	\put(76.773,24.981){\line(1,0){18}}
	\put(76.699,24.977){\line(0,-1){14.125}}
	\put(94.574,10.602){\line(0,1){14.625}}
	\qbezier(76.591,25.226)(85.946,31.743)(94.881,24.806)
	\put(76.801,10.721){\line(1,0){18.079}}
	\qbezier(76.381,10.721)(86.262,2.732)(94.46,10.721)
	\put(101.182,10.489){\circle*{1.682}}
	\put(119.051,10.456){\circle*{1.682}}
	\put(101.506,24.743){\circle*{1.682}}
	\put(119.375,24.711){\circle*{1.682}}
	\put(101.364,24.749){\line(1,0){18}}
	\put(101.29,24.745){\line(0,-1){14.125}}
	\put(119.165,10.37){\line(0,1){14.625}}
	\put(101.392,10.489){\line(1,0){18.079}}
	\put(110.845,33.517){\circle*{1.682}}
	\multiput(110.896,33.717)(.0337004049,-.0361133603){247}{\line(0,-1){.0361133603}}
	\multiput(101.27,24.854)(.0364583333,.0336174242){264}{\line(1,0){.0364583333}}
	\multiput(101.02,10.604)(.0336700337,.0778619529){297}{\line(0,1){.0778619529}}
	\put(134.05,11.086){\circle*{1.682}}
	\put(151.919,11.053){\circle*{1.682}}
	\put(134.374,25.34){\circle*{1.682}}
	\put(152.243,25.308){\circle*{1.682}}
	\put(134.232,25.346){\line(1,0){18}}
	\put(134.26,11.086){\line(1,0){18.079}}
	\multiput(134.155,25.381)(.0418962264,-.0337146226){424}{\line(1,0){.0418962264}}
	\put(126.05,19.063){\circle*{1.682}}
	\put(161.3,19.188){\circle*{1.682}}
	\multiput(152.084,11.097)(.038381743,.033713693){241}{\line(1,0){.038381743}}
	\put(161.334,19.222){\line(-3,2){9.375}}
	\multiput(125.959,19.347)(.047857143,.033571429){175}{\line(1,0){.047857143}}
	\multiput(125.834,18.972)(.035326087,-.033695652){230}{\line(1,0){.035326087}}
	\put(4.841,8.989){\circle*{1.682}}
	\put(22.71,8.956){\circle*{1.682}}
	\put(5.165,23.243){\circle*{1.682}}
	\put(23.034,23.211){\circle*{1.682}}
	\put(5.023,23.249){\line(1,0){18}}
	\put(4.949,23.245){\line(0,-1){14.125}}
	\put(22.824,8.87){\line(0,1){14.625}}
	\put(5.051,8.989){\line(1,0){18.079}}
	\put(14.504,32.017){\circle*{1.682}}
	\multiput(14.555,32.217)(.0337004049,-.0361133603){247}{\line(0,-1){.0361133603}}
	\multiput(4.929,23.354)(.0364583333,.0336174242){264}{\line(1,0){.0364583333}}
	\multiput(4.679,9.104)(.0336700337,.0778619529){297}{\line(0,1){.0778619529}}
	\multiput(14.25,32)(.0336734694,-.093877551){245}{\line(0,-1){.093877551}}
	\multiput(22.5,9)(-.0419621749,.0336879433){423}{\line(-1,0){.0419621749}}
	\multiput(4.75,9.25)(.0439759036,.0337349398){415}{\line(1,0){.0439759036}}
	\put(27.841,9.489){\circle*{1.682}}
	\put(45.71,9.457){\circle*{1.682}}
	\put(28.165,23.743){\circle*{1.682}}
	\put(46.034,23.711){\circle*{1.682}}
	\put(28.023,23.749){\line(1,0){18}}
	\put(27.949,23.745){\line(0,-1){14.125}}
	\put(45.824,9.37){\line(0,1){14.625}}
	\qbezier(27.841,23.994)(37.196,30.511)(46.131,23.574)
	\put(28.051,9.489){\line(1,0){18.079}}
	\qbezier(27.631,9.489)(37.512,1.5)(45.71,9.489)
	\multiput(27.946,23.784)(.0418962264,-.0337146226){424}{\line(1,0){.0418962264}}
	\multiput(27.631,9.489)(.043384434,.0337146226){424}{\line(1,0){.043384434}}
	\put(51.091,9.739){\circle*{1.682}}
	\put(68.96,9.707){\circle*{1.682}}
	\put(51.415,23.993){\circle*{1.682}}
	\put(69.284,23.961){\circle*{1.682}}
	\put(51.273,23.999){\line(1,0){18}}
	\put(51.199,23.995){\line(0,-1){14.125}}
	\put(69.074,9.62){\line(0,1){14.625}}
	\put(51.301,9.739){\line(1,0){18.079}}
	\multiput(51.196,24.034)(.0418962264,-.0337146226){424}{\line(1,0){.0418962264}}
	\multiput(50.881,9.739)(.043384434,.0337146226){424}{\line(1,0){.043384434}}
	\multiput(110.5,33.75)(.0336734694,-.0948979592){245}{\line(0,-1){.0948979592}}
	\multiput(125.75,19.25)(.109243697,-.033613445){238}{\line(1,0){.109243697}}
	\put(151.75,11.25){\line(0,1){14.5}}
	\multiput(134,11.25)(.117391304,.033695652){230}{\line(1,0){.117391304}}
	\qbezier(69,24.25)(63.125,29.125)(70.75,31.5)
	\qbezier(70.75,31.5)(78,27.875)(69.25,23.75)
	\qbezier(94.75,24.75)(100.875,17.75)(94.5,10.75)
	\qbezier(101,24.75)(102,34.875)(111,33.5)
	\qbezier(111,33.5)(120.25,33.375)(119.5,24.75)
	\put(13,1.5){\makebox(0,0)[cc]{$Q_1$}}
	\put(36.75,1.5){\makebox(0,0)[cc]{$Q_2$}}
	\put(60.25,1.5){\makebox(0,0)[cc]{$Q_3$}}
	\put(85.75,1.5){\makebox(0,0)[cc]{$Q_4$}}
	\put(109.75,1.5){\makebox(0,0)[cc]{$Q_5$}}
	\put(142.75,1.5){\makebox(0,0)[cc]{$Q_6$}}
\end{picture}
\caption{Forbidden Minors for the Class of Graphic Matroid which yields Graphic Matroid.}
\label{gtg}
\end{figure}

In this paper, we obtain excluded minors for the class $\mathcal{M}_3$. 
 The following is the main result of the paper.
\begin{thm}\label{mt}
	Let $M$ be a graphic matroid without containing a minor isomorphic to any member of the set $\mathcal{T}_3.$ Then $M \in \mathcal{M}_3$ if and only if $M$ does not contain a minor isomorphic to any of the matroid $\tilde{M}(G_{1})$, $\tilde{M}(G_{3})$, $M(G_4)$,  $M(G_5)$, $M(G_6)$ and $M(G_7)$,  where $G_i$ is the graph shown in Figure \ref{1} or Figure \ref{2}, for $i=1,3,4,5,6,7$.
\end{thm}

\begin{figure}[h!]
	\centering
\unitlength 1mm 
\linethickness{0.4pt}
\ifx\plotpoint\undefined\newsavebox{\plotpoint}\fi 
\begin{picture}(107.374,32.733)(0,0)
	\put(3.569,8.056){\circle*{1.682}}
	\put(21.438,8.024){\circle*{1.682}}
	\put(3.893,22.31){\circle*{1.682}}
	\put(21.762,22.278){\circle*{1.682}}
	\put(3.751,22.316){\line(1,0){18}}
	\put(3.677,22.312){\line(0,-1){14.125}}
	\put(21.552,7.937){\line(0,1){14.625}}
	\put(3.779,8.056){\line(1,0){18.079}}
	\multiput(3.674,22.351)(.0418962264,-.0337146226){424}{\line(1,0){.0418962264}}
	\multiput(3.359,8.056)(.043384434,.0337146226){424}{\line(1,0){.043384434}}
	\put(33.11,8.001){\circle*{1.682}}
	\put(50.979,7.968){\circle*{1.682}}
	\put(33.434,22.255){\circle*{1.682}}
	\put(51.303,22.223){\circle*{1.682}}
	\put(33.292,22.261){\line(1,0){18}}
	\put(33.218,22.257){\line(0,-1){14.125}}
	\put(51.093,7.882){\line(0,1){14.625}}
	\put(33.32,8.001){\line(1,0){18.079}}
	\qbezier(32.9,8.001)(42.78,.012)(50.979,8.001)
	\qbezier(33.411,22.274)(25.809,16.97)(33.057,8.132)
	\qbezier(33.411,22.45)(42.603,30.671)(51.088,22.274)
	\qbezier(21.567,22.45)(15.379,30.052)(24.749,28.461)
	\qbezier(24.749,28.461)(29.698,24.307)(21.92,22.627)
	\put(58.293,8.166){\circle*{1.682}}
	\put(76.162,8.133){\circle*{1.682}}
	\put(58.617,22.42){\circle*{1.682}}
	\put(76.486,22.388){\circle*{1.682}}
	\put(58.401,22.422){\line(0,-1){14.125}}
	\put(76.276,8.047){\line(0,1){14.625}}
	\put(58.503,8.166){\line(1,0){18.079}}
	\put(67.956,31.194){\circle*{1.682}}
	\multiput(68.007,31.394)(.0337004049,-.0361133603){247}{\line(0,-1){.0361133603}}
	\multiput(58.381,22.531)(.0364583333,.0336174242){264}{\line(1,0){.0364583333}}
	\multiput(58.131,8.281)(.0336700337,.0778619529){297}{\line(0,1){.0778619529}}
	\multiput(67.999,31.375)(.033613174,-.09716308){238}{\line(0,-1){.09716308}}
	\multiput(75.999,8.25)(-.0419618356,.0336876708){423}{\line(-1,0){.0419618356}}
	\multiput(58.125,8.125)(.0428483532,.0336876708){423}{\line(1,0){.0428483532}}
	\put(83.34,8.864){\circle*{1.682}}
	\put(101.209,8.831){\circle*{1.682}}
	\put(83.664,23.118){\circle*{1.682}}
	\put(101.533,23.086){\circle*{1.682}}
	\put(101.323,8.745){\line(0,1){14.625}}
	\put(83.55,8.864){\line(1,0){18.079}}
	\qbezier(83.13,8.864)(93.01,.875)(101.209,8.864)
	\put(93.003,31.892){\circle*{1.682}}
	\multiput(93.054,32.092)(.0337004049,-.0361133603){247}{\line(0,-1){.0361133603}}
	\multiput(83.428,23.229)(.0364583333,.0336174242){264}{\line(1,0){.0364583333}}
	\multiput(83.178,8.979)(.0336700337,.0778619529){297}{\line(0,1){.0778619529}}
	\multiput(83.499,23.125)(.0426607528,-.0337109445){419}{\line(1,0){.0426607528}}
	\qbezier(101.374,23)(107.374,17.375)(101.374,8.75)
	\qbezier(83.499,23.25)(96.374,21.875)(101.249,9)
	\put(12.193,1){\makebox(0,0)[cc]{$G_4$}}
	\put(42.255,1){\makebox(0,0)[cc]{$G_5$}}
	\put(67.271,1){\makebox(0,0)[cc]{$G_6$}}
	\put(93.759,1){\makebox(0,0)[cc]{$G_7$}}
	\put(12.613,23.755){\makebox(0,0)[cc]{}}
	\put(12.193,9.46){\makebox(0,0)[cc]{}}
	\put(24.806,30.272){\makebox(0,0)[cc]{}}
	\put(42.465,28.38){\makebox(0,0)[cc]{}}
	\put(27.96,15.767){\makebox(0,0)[cc]{}}
	\put(42.044,5.886){\makebox(0,0)[cc]{}}
	\put(61.385,27.96){\makebox(0,0)[cc]{}}
	\put(73.578,28.17){\makebox(0,0)[cc]{}}
	\put(66.851,6.727){\makebox(0,0)[cc]{}}
	\put(106.162,16.608){\makebox(0,0)[cc]{}}
	\put(94.81,20.392){\makebox(0,0)[cc]{}}
	\put(92.708,6.307){\makebox(0,0)[cc]{}}
\end{picture}

	\caption{Forbidden Minors for the Class $\mathcal{M}_3$.}
	\label{2}
\end{figure}

We prove the main result in Section 5. We also give an alternate and short proof to the characterization of the class $\mathcal{M}_2$ in the fourth section. Section 2 and Section 3 contain results that are useful to prove the paper's main result. 

\section{Minimal Matroid}
In this section, we prove that if the splitting of a graphic matroid is not cographic, then it contains a special type of graphic matroid as a minor.  Let $\mathcal{F}=\{F_{7}, F_{7}^{*}, M(K_5), M(K_{3,3})\}.$

We use the following result from Oxley's book \cite{ox}.
\begin{thm}\cite{ox}\label{cg}
	Let $M$ be a binary matroid. Then $M$ is a cographic matroid if and only if $M$ does not contain a minor isomorphic to any member of the set $\mathcal{F}$. 
\end{thm}

The following result is proved in \cite{gm}.  
\begin{lem}\label{bp1} \cite{gm}
	Let $M$ be a binary matroid and  $T \subseteq E(M)$. Then the following holds.
	\begin{enumerate}[(i).]
	\item If $x \in E(M)-T$, then $(M_T)\backslash x = (M \backslash x)_T$ and $(M_T)/x = (M/x)_T$.
	\item If $y \in T$, then  $(M_T)\backslash y= ( (M)\backslash y)_{T-{y}}$.
	\item $M_T \backslash T = M \backslash T$.
	\item $M_T = M$ , if $T$ is a cocircuit of $M$.
	\end{enumerate}
	
\end{lem}
In the following result, we prove that if the splitting matroid of a graphic matroid is not cographic, then it contains a minimal minor whose splitting is not cographic. 
\begin{lem} \label{lm1}
	Let $M$ be a graphic matroid and $k=2,3$.  If $M_{T}$ contains a minor isomorphic to $F$ for some $F \in \mathcal{F}$ and  $T \subseteq E(M)$ with $|T|=k$, then $M$ has a minor $N$ with $T \subseteq N$, such that one of the following holds.
	\begin{enumerate}[(i).]
	\item $N_{T} \cong F$.
	\item $N_{T} / T^{\prime} \cong F$ for some non-empty subset $T^{\prime}$ of $T$.
	\item $N$ is isomorphic to one element extension of forbidden minors for the class $\mathcal{M}_{k-1}$.
	\item $N$ is isomorphic to some member of $\mathcal{T}_k.$
	\end{enumerate}
	
\end{lem} 
\begin{proof}
	Let $M_T$ contains a minor isomorphic to $F$ for some $F \in \mathcal{F}$ and $T \subseteq E(M)$ with $|T|=k$, where $k=2,3$. Then there exist $T_1, T_2 \subseteq E(M)$ such that $M_T \backslash T_1 / T_2 \cong F$.	Now let $T_i^{\prime}=T \cap T_i$  and $T_i^{\prime\prime}=T_i - T_i^{\prime}$ for $i=1,2$. Then $T_i^{\prime} \subseteq T$ and $T_i^{\prime\prime} \cap T= \phi$.  By Lemma \ref{bp1}(1), $ M_T\backslash T_1^{\prime\prime}/T_2^{\prime\prime} = (M\backslash T_1^{\prime\prime}/T_2^{\prime\prime})_T$. Let $N=M\backslash T_1^{\prime\prime}/T_2^{\prime\prime}$. Then $N$ is minor of $M$ containing $T$ and  $N_T\backslash T_1^{\prime}/T_2^{\prime}=(M\backslash T_1^{\prime\prime}/T_2^{\prime\prime})_T \backslash T_1^{\prime}/T_2^{\prime}= M_T\backslash T_1^{\prime\prime}/T_2^{\prime\prime}\backslash T_1^{\prime}/T_2^{\prime}=M_T\backslash T_1/T_2\cong F.$

	If $T_1^{\prime}=\phi $ and $T_2^{\prime}=\phi,$ then (i) or (iv)  holds.   If $T_1^{\prime}=\phi $ and $T_2^{\prime} \neq \phi$ then (ii) or (iv) holds. If $T_1^{\prime} \neq \phi $, then  $|T_1^{\prime}|\leq k=|T|$.
	Suppose $|T_1^{\prime}|=k.$ Then  $T_1^{\prime}=T$ and $F\cong N_T \backslash T_1^{\prime}=N_T \backslash T=N\backslash T$ this shows that $F$ is minor of graphic matroid $M$ which shows that (iv) holds. Let $|T_1^{\prime}|< k$ and $z \in T_1^{\prime}.$ Suppose $N\backslash z=P$. 
Then, by Lemma \ref{lm1} (ii), 	$N_T \backslash T_1^{\prime} /T_2^{\prime}=N_T\backslash z \backslash T_1^{\prime}-\{z\}  \backslash T_2^{\prime}=P_{T-\{z\}} \backslash T_1^{\prime}-\{z\} /T_2^{\prime} \cong F.$
	 Therefore $N$ is extension of $P,$ which is a forbidden minor for the class  $\mathcal{M}_{k-1}.$ Which means (iii) hold.  
\end{proof}

We focus on finding forbidden minors for the class $\mathcal{M}_k$ other than the members of  $\mathcal{T}_k$ for $k=2,3.$ Henceforth, we consider the graphic matroids, which do not contain a minor isomorphic to any member of $\mathcal{T}_k$ for $k=2,3.$

\begin{lem}\label{lm2}
	Let $N$ be the minor of the matroid $M$ stated in Lemma \ref{lm1} and without containing a minor isomorphic to any member of $\mathcal{T}_k$. Then $N$ does not contain 2-cocircuit. 
\end{lem}
\begin{proof}
	The proof is similar to the proof of Lemma 3.3.2, in \cite{gm}.
\end{proof}

\begin{lem}\label{lm3} 
	Let $N$ be the minor of the matroid $M$ stated in Lemma \ref{lm1}, such that $N$ does not satisfy conditions (iii) and (iv). Then $N$ does not contain a coloop.
\end{lem}
\begin{proof}
	The proof is similar to the proof of Lemma 3.3.3, in \cite{gm}.
\end{proof}

We call the matroid $N$ of Lemmas \ref{lm1}, \ref{lm2} and \ref{lm3} as a $k$-minimal matroid with repsect to $F$ where $F \in \mathcal{F}.$ We define a $k$-minimal matroid with respect $F \in \mathcal{F}$ as follows.
\begin{defn}\label{mini}
	Let $k=2,3$ and $N$ be a graphic matroid without containing any member of $\mathcal{T}_k. $  Then $N$ is called $k$-minimal matroid with respect to some $F$ of $\mathcal{F}$ if the following conditions hold.
	\begin{enumerate}[(i).]
		\item $N$ contains no coloop and no 2-cocircuit.
		\item $N$ is not a single element extension of any member of forbidden minors of $\mathcal{M}_{k-1}$. 
		\item There is a set $T \subseteq E(N)$ with $|T|=k$  such that $N_{T} \cong F$ or $N_{T} / T^{\prime} \cong F$ for some non-empty subset $T^{\prime}$ of $T$.
	\end{enumerate}
	
\end{defn}
In light of Lemma \ref{lm1} and Definition \ref{mini}, we focus on to find the $k$-minimal matroids with respect to every member of $\mathcal{F}$, for $k=2,3$.
Let $M$ be a  matroid. Then matroid $Q$ is called a quotient of $M$ if there exist a matroid $N$ such that $Q=N/a$ and $N\backslash a=M$ for some $a \in N$. Here $M$ is called a lift of $Q$. We prove that a $k$-minimal matroid with respect to $F$ contains a minor isomorphic to a quotient of $F$ in the following lemma. 

\begin{lem}\label{ql}
	Let $k=2,3$ and let $M$ be a $k$-minimal matroid  with respect to $F \in \mathcal{F}$. Then there is a binary matroid $N$ containing an element $a$ such that $N \backslash a \cong F$. Further, $N/a$ is a minor of $M$ such that either $M=N/a$ or $M$ is a coextension of $N/a$ by at most $k$ elements.
\end{lem}

\begin{proof}
	Since $M$ is a $k$-minimal matroid, $M$ is graphic and $M_T \cong F$ or $M_T/T' \cong F$ for some $T \subseteq E(M)$ and $T' \subseteq T$ with $|T|=k$. Let $A$ be the standard matrix representation of $M$ over the field $GF(2)$. Let $A_{T}^{\prime}$ be the matrix obtained from $A$ by adding one extra row to the matrix $A$ at the bottom whose entries are $1$ in the columns corresponding to the elements of $T$ and zero otherwise and then adding one extra column labeled $a$ with entry $1$ in the last row and $0$ elsewhere. Let $M_{T}^{\prime}$ be the vector matroid of $A_{T}^{\prime}$. Here $E(M_{T}^{\prime})= E(M) \cup\{a\}$. Therefore 	(i). $M_{T}=M_{T}^{\prime} \backslash a$ and
	(ii). $M_{T}^{\prime} / a \cong M.$ 
	
	\vskip.3cm

	\begin{figure}[h!]
	\centering
	\unitlength 1mm 
	\linethickness{0.4pt}
	\ifx\plotpoint\undefined\newsavebox{\plotpoint}\fi 
	\begin{picture}(33.5,26.25)(0,0)
		\put(30.75,26.25){\makebox(0,0)[cc]{$N=M_T'$}}
		\put(30.25,4){\makebox(0,0)[cc]{$M$}}
		\put(10.5,16.25){\makebox(0,0)[cc]{$M_T$}}
		\put(29.75,24.5){\vector(0,-1){16}}
		\put(14.75,17.75){\vector(-3,-2){.07}}\multiput(26.25,25)(-.053488372,-.03372093){215}{\line(-1,0){.053488372}}
		\put(13,14.25){\vector(-3,2){.07}}\multiput(27.75,5.25)(-.0552434457,.0337078652){267}{\line(-1,0){.0552434457}}
		\put(17.75,24.5){\makebox(0,0)[cc]{$N\backslash a$}}
		\put(33.5,17.25){\makebox(0,0)[cc]{$N/a$}}
	\end{picture}
	\end{figure}
	
	Suppose $M_T \cong F.$  Let $N=M_T'.$ Then $N$ is a binary matroid containing $a$, such that $M_T'/a \cong N/a \cong M$ and $M_T'\backslash a \cong M_T \cong F \cong N \backslash a $.

	Suppose $M_T/T' \cong F$ for some $T' \subset T.$ Let $N=M_T'/T'$. Then $N \backslash a \cong M_T'/T' \backslash a \cong M_T'\backslash a/T'\cong M_T/T' \cong F$ and $N/a \cong M_T'/T'/a \cong M_T'/a/T' \cong M/T'.$ Here $M$ is a coextension of $N/a$ by at most $k$ elements as $|T'| \leq |T|= k$.
\end{proof}

In the above result, $N/a$ is a quotient of $F$, and $M$ is isomorphic to a quotient of $F$ or a coextension of a quotient of $F$. Note that $N$ is binary, and we obtain quotients through binary matroid coextension/extension. Since $M$ is graphic, in the next section, we determine graphic quotients of $F \in \mathcal{F}$.

\section{Graphic Quotients of Non-Cographic Matroids}

In this section, we obtain graphic quotients of every member of the set $\mathcal{F}.$  First, we prove the following lemma. 

\begin{lem}\label{glem}
	Let $N$ be a binary matroid containing an element $a$ such that $N \backslash a \cong F$ for some $F \in \mathcal{F}$ and $N/a \cong M(G)$ for some connected graph $G$. Then the following holds.
	\begin{enumerate}[(i).]
		\item The graph $G$ is a block, or it has at most two blocks, one of which is a loop.
		\item The graph $G$ contains at most one loop.
		\item The graph $G$ does not contain more than two parallel edges with the same end vertices.
		\item The graph $G$ does not contain a 2-edge cut.
		\item If $F$ is Eulerian, then G is also Eulerian. 
	\end{enumerate}
\end{lem}
\begin{proof}
	The proof is similar to the proof of Lemma 3.4.1, in \cite{gm}.
\end{proof}

Mundhe \cite{gm}, determined graphic quotients for $F_7^*$ and $F_7$ as follows.
\begin{lem}\cite{gm} \label{q1}
Every graphic quotient of $F_7$ is isomorphic to the circuit matroid $M(H_3)$ and every graphic quotient of $F_7^*$ is isomorphic to the circuit matroid $M(H_1)$ or $M(H_2),$ where $H_1$, $H_2$, $H_3$ are shown in Figure \ref{3}.
\end{lem}
\begin{figure}[h!]
\centering		
\unitlength 1mm 
\linethickness{0.4pt}
\ifx\plotpoint\undefined\newsavebox{\plotpoint}\fi 
\begin{picture}(79.358,32.353)(0,0)
	\put(3.359,9.738){\circle*{1.682}}
	\put(21.228,9.706){\circle*{1.682}}
	\put(3.683,23.992){\circle*{1.682}}
	\put(21.552,23.96){\circle*{1.682}}
	\put(3.541,23.998){\line(1,0){18}}
	\put(3.467,23.994){\line(0,-1){14.125}}
	\put(21.342,9.619){\line(0,1){14.625}}
	\put(3.569,9.738){\line(1,0){18.079}}
	\multiput(3.464,24.033)(.0418962264,-.0337146226){424}{\line(1,0){.0418962264}}
	\multiput(3.149,9.738)(.043384434,.0337146226){424}{\line(1,0){.043384434}}
	\put(32.9,9.683){\circle*{1.682}}
	\put(50.769,9.65){\circle*{1.682}}
	\put(33.224,23.937){\circle*{1.682}}
	\put(51.093,23.905){\circle*{1.682}}
	\put(33.082,23.943){\line(1,0){18}}
	\put(33.008,23.939){\line(0,-1){14.125}}
	\put(50.883,9.564){\line(0,1){14.625}}
	\put(33.11,9.683){\line(1,0){18.079}}
	\qbezier(32.69,9.683)(42.57,1.694)(50.769,9.683)
	\qbezier(33.201,23.956)(25.599,18.652)(32.847,9.814)
	\qbezier(33.201,24.132)(42.393,32.353)(50.878,23.956)
	\qbezier(21.357,24.132)(15.169,31.734)(24.539,30.143)
	\qbezier(24.539,30.143)(29.488,25.989)(21.71,24.309)
	\put(58.862,10.543){\circle*{1.682}}
	\put(76.731,10.511){\circle*{1.682}}
	\put(59.072,10.543){\line(1,0){18.079}}
	\put(66.85,25.437){\circle*{1.682}}
	\multiput(58.652,10.301)(.033739164,.063152793){243}{\line(0,1){.063152793}}
	\multiput(66.85,25.647)(.0337214481,-.0502234333){293}{\line(0,-1){.0502234333}}
	\qbezier(58.652,10.721)(55.078,23.545)(66.64,25.437)
	\qbezier(76.731,10.721)(64.433,.736)(58.441,10.511)
	\qbezier(66.64,25.647)(79.358,22.914)(76.52,10.931)
	\qbezier(66.64,25.437)(58.231,30.062)(66.64,32.164)
	\qbezier(66.64,32.164)(72.842,30.377)(66.85,25.647)
	\put(10.511,.841){\makebox(0,0)[cc]{$H_1$}}
	\put(42.465,.841){\makebox(0,0)[cc]{$H_2$}}
	\put(67.691,.841){\makebox(0,0)[cc]{$H_3$}}
\end{picture}

		\caption{Quotients of $F_7^*$ and $F_7$.}
		\label{3}
	\end{figure}
We find all graphic quotients of $M(K_5)$ and $M(K_{3,3}).$
\begin{lem} \label{q3}
 Every graphic quotient of $M(K_5)$ obtained via binary extension is isomorphic to the circuit matroid $M(H_{4})$, $M(H_{5})$ or $M(H_{6})$,  where $H_4$, $H_5$ and $H_{6}$ are the graphs as shown in Figure \ref{4}.
\end{lem}

\begin{figure}[h!]
	\centering
\unitlength 1mm 
\linethickness{0.4pt}
\ifx\plotpoint\undefined\newsavebox{\plotpoint}\fi 
\begin{picture}(91.026,32.731)(0,0)
	\put(31.949,9.738){\circle*{1.682}}
	\put(49.818,9.706){\circle*{1.682}}
	\put(32.273,23.992){\circle*{1.682}}
	\put(50.142,23.96){\circle*{1.682}}
	\put(32.131,23.998){\line(1,0){18}}
	\put(32.057,23.994){\line(0,-1){14.125}}
	\put(49.932,9.619){\line(0,1){14.625}}
	\put(32.159,9.738){\line(1,0){18.079}}
	\multiput(32.054,24.033)(.0418962264,-.0337146226){424}{\line(1,0){.0418962264}}
	\multiput(31.739,9.738)(.043384434,.0337146226){424}{\line(1,0){.043384434}}
	\qbezier(32.164,24.175)(39.206,32.269)(50.033,23.965)
	\qbezier(50.033,23.965)(58.126,17.133)(49.822,9.88)
	\qbezier(31.533,9.67)(45.933,12.403)(49.822,23.545)
	\put(64.328,10.123){\circle*{1.682}}
	\put(82.197,10.091){\circle*{1.682}}
	\put(64.652,24.377){\circle*{1.682}}
	\put(82.521,24.345){\circle*{1.682}}
	\put(64.51,24.383){\line(1,0){18}}
	\put(64.436,24.379){\line(0,-1){14.125}}
	\put(82.311,10.004){\line(0,1){14.625}}
	\put(64.538,10.123){\line(1,0){18.079}}
	\multiput(64.433,24.418)(.0418962264,-.0337146226){424}{\line(1,0){.0418962264}}
	\multiput(64.118,10.123)(.043384434,.0337146226){424}{\line(1,0){.043384434}}
	\qbezier(64.543,24.56)(71.585,32.654)(82.412,24.35)
	\qbezier(82.412,24.35)(90.505,17.518)(82.201,10.265)
	\qbezier(63.912,10.055)(78.312,12.788)(82.201,23.93)
	\qbezier(49.612,9.67)(56.444,14.4)(55.288,6.096)
	\qbezier(55.288,6.096)(50.033,.21)(49.822,9.46)
	\qbezier(82.617,24.386)(80.304,31.428)(87.242,29.641)
	\qbezier(87.242,29.641)(91.026,23.755)(83.037,24.596)
	\put(4.835,8.862){\circle*{1.682}}
	\put(22.704,8.829){\circle*{1.682}}
	\put(5.159,23.116){\circle*{1.682}}
	\put(23.028,23.084){\circle*{1.682}}
	\put(4.943,23.118){\line(0,-1){14.125}}
	\put(22.818,8.743){\line(0,1){14.625}}
	\put(14.498,31.89){\circle*{1.682}}
	\multiput(14.549,32.09)(.0337004049,-.0361133603){247}{\line(0,-1){.0361133603}}
	\multiput(4.923,23.227)(.0364583333,.0336174242){264}{\line(1,0){.0364583333}}
	\multiput(4.784,23.286)(.0420257009,-.0336892523){428}{\line(1,0){.0420257009}}
	\multiput(4.784,9.164)(.0436987952,.0336698795){415}{\line(1,0){.0436987952}}
	\put(4.625,23.335){\line(1,0){18.289}}
	\put(4.625,8.619){\line(1,0){18.079}}
	\multiput(22.704,8.619)(-.0336353508,.0941789821){250}{\line(0,1){.0941789821}}
	\multiput(14.295,32.164)(-.0336652754,-.08378913){281}{\line(0,-1){.08378913}}
	\put(13.664,3.994){\makebox(0,0)[cc]{$H_4$}}
	\put(40.783,3.994){\makebox(0,0)[cc]{$H_5$}}
	\put(73.788,3.994){\makebox(0,0)[cc]{$H_6$}}
\end{picture}

	\caption{Quotients of $M(K_5)$}
	\label{4}
\end{figure}
\begin{proof} 
	Let $N$ be binary matroid containing $a$ such that $N \backslash a \cong M(K_5)$ and $N/a$ is graphic. Assume that $N/a=M(G)$ for a connected graph $G$. Suppose $a$ is a loop or a coloop of $N$.  Then $N/a = N\backslash a \cong M(K_5)$. Thus $G \cong H_4$.
	
	Suppose $a$ is not a loop or a coloop of $N$. Since the graph $K_5$ has 5 vertices and 10 edges, $r(N)=4$ and $|E(N)|=11$. Hence $r(N/a)=3$ and $|E(N/a)|=10.$  Thus graph $G$ corresponding to $N/a$ has 4 vertices and 10 edges. By Lemma \ref{glem}(v), $G$ is Eulerian as $K_5$ is Eulerian and by Lemma \ref{glem}, $G$ contains at most one loop.

	 If $G$ does not contain a loop, then, by Harary (\cite{Har}, pg 230), there is only one Eulerian graph on $4$ vertices and $10$ edges as shown in Figure \ref{20}. It contains five distinct 2-cycles. That is, $N/a$ contains five 2-circuits. One can easily check that $N\backslash a$ contains an odd cocircuit or a 2-circuit. Since $N\backslash a \cong M(K_5),$  $M(K_5)$ contains an odd cocircuit or a 2-circuit, a contradiction to $M(K_5)$ is  a Eulerian graph.  Hence we discard this graph.

	\begin{figure}[h!]
		\centering
	\unitlength 1mm 
	\linethickness{0.4pt}
	\ifx\plotpoint\undefined\newsavebox{\plotpoint}\fi 
	\begin{picture}(33,27.875)(0,0)
		\put(7.207,17.332){\circle*{1.414}}
		\put(24.707,17.332){\circle*{1.414}}
		\put(24.707,6.332){\circle*{1.414}}
		\put(7.707,6.332){\circle*{1.414}}
		\thicklines
		\put(7,17.375){\line(1,0){17.75}}
		\put(24.75,17.375){\line(0,1){0}}
		\put(24.75,17.375){\line(0,-1){11.25}}
		\put(24.75,6.125){\line(0,1){.25}}
		\put(24.75,6.375){\line(-1,0){17.25}}
		\put(7.5,6.375){\line(0,1){0}}
		\multiput(7.5,6.375)(-.0333333,.7166667){15}{\line(0,1){.7166667}}
		\put(7,17.125){\line(0,1){0}}
		\qbezier(7,17.125)(16.5,27.875)(25,17.625)
		\qbezier(7.25,17.625)(-2,11.5)(7.75,5.875)
		\qbezier(8,5.875)(16.75,-3.5)(24.5,6.625)
		\thinlines
		\qbezier(25,6.5)(33,11)(25,17.5)
		\qbezier(6.75,17.75)(23.375,13.75)(24.5,6.75)
		\qbezier(7,17.5)(13.75,9.375)(24.5,6.75)
	\end{picture}

	\caption{Eulerian Graph on 4 Vertices and 10 Edges}
	\label{20}
	\end{figure}
	 If $G$ contains a loop, then $G$ minus that loop is an Eulerian graph of $4$ vertices and $9$ edges. By Harary (\cite{Har}, pg 230), there is only one Eulerian graph of $4$ vertices and $9$ edges. That graph is $H_5$ minus the loop as shown in Figure  \ref{4}. Thus $G$ can be obtained by adding a loop to that graph. Hence $G$ is isomorphic to either the graph $H_5$ or $H_6$ as shown in Figure \ref{4}.   
\end{proof}

We need the following result to find graphic quotients of $M(K_{3,3})$.

\begin{lem}\label{qlemma}
Let $N$ be a binary matroid and $a \in E(N)$. If $N\backslash a \cong M(K_{3,3})$, then no element of $N/a$ belongs to 2-circuit and odd circuit.
\end{lem} 
\begin{proof} Since $M(K_{3,3})$ is a bipartite graph, it contain only circuits with even length. On contrary, let $x_1$ be an element of  $N/a$ belongs to both 2-circuit and odd circuit. Suppose $\{x_1,x\}$ is a 2-circuit and $C_1=\{x_1, y_1, y_2, \cdots y_k\}$ is an odd circuit, where $k$ is an even integer.   Then $C_2=\{x, y_1, y_2, \cdots y_k\}$ is also a circuit of  $N/a$. Then $D_1=\{x_1, y_1, y_2, \cdots y_k,a\}$ and $D_2=\{x, y_1, y_2, \cdots y_k,a\}$ are circuits of $N$ otherwise $N\backslash a$ contains odd circuit, a contradiction to $N\backslash a$ is bipartite.  Suppose $D_1$ and $D_2$ are circuits of $N$. Then $D_1 \triangle D_2= \{x_1, x\} $ is a circuit in $N$ and also in $N\backslash a \cong M(K_{3,3})$, a contradiction. Thus no element of $N/a$ lies in a 2-circuit as well as in an odd circuit.  
\end{proof}
In this paper, we aim to find the forbidden minors for the class $\mathcal{M}_3$ and, we also intend to give an alternative proof for the   forbidden minor characterization of the class $\mathcal{M}_2$. By Theorem \ref{gct}, $M(G_2)$ is a forbidden minor for the class $\mathcal{M}_2.$ Therefore, by Lemma \ref{lm1}, $\tilde{M}({G_2})$ is a forbidden minor for the class $\mathcal{M}_3$, where $\tilde{M}({G_2})$ is an extension of $G_2$ by single element.  Thus it is sufficient, if we consider the graphic quotients  $M(K_{3,3})$ which avoids  ${\tilde{M}({G_2})}$ as a minor. 

We now find quotient for $M(K_{3,3})$.
\begin{lem} \label{q4} Every graphic quotient of $M(K_{3,3})$ obtained via binary extension and without containing a minor isomorphic to $\tilde{M}({G_2})$ is isomorphic to the circuit matroid of $H_{i}$ for some  $i \in \{7,8,9,10,11\}$, where $H_{i}$ is the graph as shown in Figure \ref{5}, for $i=7,8,9,10,11$.
\end{lem}
\begin{figure}
\centering
\unitlength 1mm 
\linethickness{0.4pt}
\ifx\plotpoint\undefined\newsavebox{\plotpoint}\fi 
\begin{picture}(142.691,37.433)(0,0)
	\put(93.61,8.586){\circle*{1.682}}
	\put(111.479,8.553){\circle*{1.682}}
	\put(93.934,22.84){\circle*{1.682}}
	\put(111.803,22.808){\circle*{1.682}}
	\put(93.718,22.842){\line(0,-1){14.125}}
	\put(111.593,8.467){\line(0,1){14.625}}
	\put(93.82,8.586){\line(1,0){18.079}}
	\put(103.273,31.614){\circle*{1.682}}
	\multiput(103.324,31.814)(.0337004049,-.0361133603){247}{\line(0,-1){.0361133603}}
	\multiput(93.698,22.951)(.0364583333,.0336174242){264}{\line(1,0){.0364583333}}
	\multiput(93.448,8.701)(.0336700337,.0778619529){297}{\line(0,1){.0778619529}}
	\multiput(103.316,31.795)(.033613445,-.097163866){238}{\line(0,-1){.097163866}}
	\put(118.657,9.284){\circle*{1.682}}
	\put(136.526,9.251){\circle*{1.682}}
	\put(118.981,23.538){\circle*{1.682}}
	\put(136.85,23.506){\circle*{1.682}}
	\put(136.64,9.165){\line(0,1){14.625}}
	\put(118.867,9.284){\line(1,0){18.079}}
	\qbezier(118.447,9.284)(128.327,1.295)(136.526,9.284)
	\put(128.32,32.312){\circle*{1.682}}
	\multiput(128.371,32.512)(.0337004049,-.0361133603){247}{\line(0,-1){.0361133603}}
	\multiput(118.745,23.649)(.0364583333,.0336174242){264}{\line(1,0){.0364583333}}
	\multiput(118.495,9.399)(.0336700337,.0778619529){297}{\line(0,1){.0778619529}}
	\multiput(118.816,23.545)(.0426610979,-.0337112172){419}{\line(1,0){.0426610979}}
	\qbezier(136.691,23.42)(142.691,17.795)(136.691,9.17)
	\qbezier(118.816,23.67)(131.691,22.295)(136.566,9.42)
	\put(41.165,8.591){\circle*{1.682}}
	\put(59.034,8.558){\circle*{1.682}}
	\put(41.489,22.845){\circle*{1.682}}
	\put(59.358,22.813){\circle*{1.682}}
	\put(41.273,22.847){\line(0,-1){14.125}}
	\put(59.148,8.472){\line(0,1){14.625}}
	\put(41.375,8.591){\line(1,0){18.079}}
	\put(50.828,31.619){\circle*{1.682}}
	\multiput(50.879,31.819)(.0337004049,-.0361133603){247}{\line(0,-1){.0361133603}}
	\multiput(41.253,22.956)(.0364583333,.0336174242){264}{\line(1,0){.0364583333}}
	\multiput(41.003,8.706)(.0336700337,.0778619529){297}{\line(0,1){.0778619529}}
	\multiput(50.871,31.8)(.033613445,-.097163866){238}{\line(0,-1){.097163866}}
	\multiput(58.871,8.675)(-.0419598109,.0336879433){423}{\line(-1,0){.0419598109}}
	\multiput(40.997,8.55)(.0428463357,.0336879433){423}{\line(1,0){.0428463357}}
	\put(67.362,8.484){\circle*{1.682}}
	\put(85.231,8.451){\circle*{1.682}}
	\put(67.686,22.738){\circle*{1.682}}
	\put(85.555,22.706){\circle*{1.682}}
	\put(67.47,22.74){\line(0,-1){14.125}}
	\put(85.345,8.365){\line(0,1){14.625}}
	\put(67.572,8.484){\line(1,0){18.079}}
	\put(77.025,31.512){\circle*{1.682}}
	\multiput(77.076,31.712)(.0337004049,-.0361133603){247}{\line(0,-1){.0361133603}}
	\multiput(67.45,22.849)(.0364583333,.0336174242){264}{\line(1,0){.0364583333}}
	\multiput(67.2,8.599)(.0336700337,.0778619529){297}{\line(0,1){.0778619529}}
	\multiput(77.068,31.693)(.033613445,-.097163866){238}{\line(0,-1){.097163866}}
	\put(67.425,22.866){\line(1,0){18.432}}
	\put(93.884,22.717){\line(1,0){18.135}}
	\qbezier(76.79,31.785)(71.067,35.501)(76.938,37.433)
	\qbezier(76.938,37.433)(81.918,35.575)(77.087,31.636)
	\qbezier(93.884,22.866)(87.864,20.859)(89.871,26.879)
	\qbezier(89.871,26.879)(94.702,28.737)(93.884,23.163)
	\put(4.415,9.493){\circle*{1.682}}
	\put(17.869,9.25){\circle*{1.682}}
	\put(4.739,23.747){\circle*{1.682}}
	\put(18.193,23.505){\circle*{1.682}}
	\put(4.523,23.749){\line(0,-1){14.125}}
	\put(17.983,9.164){\line(0,1){14.625}}
	\put(31.323,9.04){\circle*{1.682}}
	\put(31.647,23.295){\circle*{1.682}}
	\put(31.437,8.954){\line(0,1){14.625}}
	\multiput(4.415,9.46)(.0337390256,.0347771495){405}{\line(0,1){.0347771495}}
	\multiput(18.079,23.545)(.0336786696,-.0363946913){387}{\line(0,-1){.0363946913}}
	\multiput(31.113,9.46)(-.062088257,.03373299){430}{\line(-1,0){.062088257}}
	\multiput(4.415,23.965)(.0336994026,-.0363739583){393}{\line(0,-1){.0363739583}}
	\multiput(17.658,9.67)(.0342580876,.0337390256){405}{\line(1,0){.0342580876}}
	\put(3.994,9.67){\line(2,1){27.749}}
	\put(18.079,1){\makebox(0,0)[cc]{$H_7$}}
	\put(50.032,1){\makebox(0,0)[cc]{$H_8$}}
	\put(76.1,1){\makebox(0,0)[cc]{$H_9$}}
	\put(103.218,1){\makebox(0,0)[cc]{$H_{10}$}}
	\put(127.604,1){\makebox(0,0)[cc]{$H_{11}$}}
\end{picture}

\caption{Quotients of $M(K_{3,3})$}
	\label{5}
	\end{figure}
\begin{proof} Let $N$ be a binary matroid and $a \in E(N)$ such that $N\backslash a \cong M(K_{3,3})$ and $N/a$ is graphic.  Suppose $N/a=M(G)$ for some connected graph $G$. Assume that $M(G)$ does not contain a $\tilde{M}(G_2)$ as a minor.	If $a$ is a loop or coloop of $N,$ then $N\backslash a=N/a=K_{3,3}$. Thus $G \cong H_7$.

Suppose $a$ is not a loop or a coloop. Since the graph $K_{3,3}$ has 6 vertices and 9 edges, $r(N)=5$ and $|E(N)|=10$. Thus $r(N/a)=4$ and $|E(N/a)|=9$. Hence $N/a$ corresponds to a graph $G$ on 5 vertices and 9 edges.  By Lemma \ref{glem}(1) and (4), $G$ contains at most two block one of which is a loop and  $G$ does not contain 2-edge cut. Suppose $G$ is  a simple graph. As $G$ is 2-connected, by Harary (\cite{Har}, Pg. 217), $H_8$ is only 2-connected simple graph as shown in Figure $\ref{5}$. Thus $G \cong H_8$.

 Suppose $G$ is not a simple graph. Then $G$ contains a loop or 2-cycle. By Lemma \ref{glem}[1], $G$ is 2-connected, or $G$ can be obtained by adding one loop or parallel edges to 2-connected or both. Thus there are the following cases.
(i). $G$ contains either a loop or a 2-cycle. (ii). $G$ contains either a loop and a 2-cycle or two 2-cycles. (iii). $G$ contains a loop and two 2-cycles or three 2-cycle. (iv) $G$ contains a loop and three 2-cycles or four 2-cycles. (v) $G$ contains a loop and four 2-cycles or at least five 2-cycles. We consider all these cases one by one.

	{\bf Case (i).} Suppose $G$ contains either one loop or one  2-cycle. 

\noindent Then $G$ can be obtained by adding a loop or a parallel edge to a 2-connected simple graph on $5$ vertices and $8$ edges. By Harary (\cite{Har},Pg 217), there are two non-isomorphic simple graphs $A_i$ and $A_{ii}$ on $5$ vertices and $8$ as shown in   Figure \ref{6}.
		
\begin{figure}[h!]
\centering
\unitlength 1mm 
\linethickness{0.4pt}
\ifx\plotpoint\undefined\newsavebox{\plotpoint}\fi 
\begin{picture}(51.079,32.04)(0,0)
	\put(5.848,8.171){\circle*{1.682}}
	\put(23.717,8.138){\circle*{1.682}}
	\put(6.172,22.425){\circle*{1.682}}
	\put(24.041,22.393){\circle*{1.682}}
	\put(5.956,22.427){\line(0,-1){14.125}}
	\put(23.831,8.052){\line(0,1){14.625}}
	\put(6.058,8.171){\line(1,0){18.079}}
	\put(15.511,31.199){\circle*{1.682}}
	\multiput(15.562,31.399)(.0337004049,-.0361133603){247}{\line(0,-1){.0361133603}}
	\multiput(5.936,22.536)(.0364583333,.0336174242){264}{\line(1,0){.0364583333}}
	\multiput(23.554,8.255)(-.0419598109,.0336879433){423}{\line(-1,0){.0419598109}}
	\multiput(5.68,8.13)(.0428463357,.0336879433){423}{\line(1,0){.0428463357}}
	\put(32.045,8.064){\circle*{1.682}}
	\put(49.914,8.031){\circle*{1.682}}
	\put(32.369,22.318){\circle*{1.682}}
	\put(50.238,22.286){\circle*{1.682}}
	\put(32.153,22.32){\line(0,-1){14.125}}
	\put(50.028,7.945){\line(0,1){14.625}}
	\put(32.255,8.064){\line(1,0){18.079}}
	\put(41.708,31.092){\circle*{1.682}}
	\multiput(41.759,31.292)(.0337004049,-.0361133603){247}{\line(0,-1){.0361133603}}
	\multiput(32.133,22.429)(.0364583333,.0336174242){264}{\line(1,0){.0364583333}}
	\multiput(31.883,8.179)(.0336700337,.0778619529){297}{\line(0,1){.0778619529}}
	\multiput(41.751,31.273)(.033613445,-.097163866){238}{\line(0,-1){.097163866}}
	\put(5.946,22.595){\line(1,0){18.135}}
	\put(32.108,22.446){\line(1,0){18.284}}
	\put(14.568,4.311){\makebox(0,0)[cc]{$A_i$}}
	\put(40.581,4.311){\makebox(0,0)[cc]{$A_{ii}$}}
\end{picture}

			\caption{Non-Isomorphic Simple Graphs on $5$ vertices and $8$ edges}
			\label{6}
\end{figure}
	Suppose $G$ is obtained from $A_i$ by adding a loop.  Then $G$ contains a  2-edge cut, by Lemma \ref{glem} (4), which is a contradiction. Hence $G$ is obtained from $A_i$ by adding a 2-cycle such that $G$ does not contain a 2-edge cut. Then $M(G) \cong {M(\tilde{G_2})},$ a contradiction. Hence we discard $A_i$.

	If $G$ is obtained from $A_{ii}$ by adding a loop, then $G \cong H_9$ or $G \cong H_{10}$ shown in Figure \ref{5}. Suppose $G$ is obtained from $A_{ii}$ by adding a edge parallel to some edge of $A_{ii}.$ Then either $M(G) \cong \tilde{M}({G_2})$ or by Lemma \ref{qlemma}, $N \backslash a$ contains a $3$-circuit or 2-circuit, which is a contradiction.  Hence we discard it.

	{\bf Case (ii).} Suppose $G$ contains either a loop and a 2-cycle or two 2-cycles. 
	
	\noindent $G$ can be obtained from a 2-connected simple graph on $5$ vertices and $7$ edges by adding a loop and a parallel edge or $2$ parallel edges. By Harary (\cite{Har}, Pg 217), there are three non-isomorphic 2-connected simple graphs $A_{iii}$, $A_{iv}$, and $A_{v}$ on 5 vertices and 7 edges as shown in Figure \ref{7}.
	\begin{center}
		\begin{figure} [h]
		\unitlength 1mm 
		\linethickness{0.4pt}
		\ifx\plotpoint\undefined\newsavebox{\plotpoint}\fi 
		\begin{picture}(76.659,32.175)(0,0)
			\put(5.848,8.171){\circle*{1.682}}
			\put(23.717,8.138){\circle*{1.682}}
			\put(6.172,22.425){\circle*{1.682}}
			\put(24.041,22.393){\circle*{1.682}}
			\put(5.956,22.427){\line(0,-1){14.125}}
			\put(23.831,8.052){\line(0,1){14.625}}
			\put(6.058,8.171){\line(1,0){18.079}}
			\put(15.511,31.199){\circle*{1.682}}
			\multiput(15.562,31.399)(.0337004049,-.0361133603){247}{\line(0,-1){.0361133603}}
			\multiput(5.936,22.536)(.0364583333,.0336174242){264}{\line(1,0){.0364583333}}
			\multiput(5.68,8.13)(.0428463357,.0336879433){423}{\line(1,0){.0428463357}}
			\put(32.045,8.064){\circle*{1.682}}
			\put(49.914,8.031){\circle*{1.682}}
			\put(32.369,22.318){\circle*{1.682}}
			\put(50.238,22.286){\circle*{1.682}}
			\put(32.153,22.32){\line(0,-1){14.125}}
			\put(50.028,7.945){\line(0,1){14.625}}
			\put(32.255,8.064){\line(1,0){18.079}}
			\put(41.708,31.092){\circle*{1.682}}
			\multiput(41.759,31.292)(.0337004049,-.0361133603){247}{\line(0,-1){.0361133603}}
			\multiput(32.133,22.429)(.0364583333,.0336174242){264}{\line(1,0){.0364583333}}
			\multiput(41.751,31.273)(.033613445,-.097163866){238}{\line(0,-1){.097163866}}
			\put(5.946,22.595){\line(1,0){18.135}}
			\put(32.108,22.446){\line(1,0){18.284}}
			\put(57.625,8.306){\circle*{1.682}}
			\put(75.494,8.273){\circle*{1.682}}
			\put(57.949,22.56){\circle*{1.682}}
			\put(75.818,22.528){\circle*{1.682}}
			\put(57.733,22.562){\line(0,-1){14.125}}
			\put(75.608,8.187){\line(0,1){14.625}}
			\put(67.288,31.334){\circle*{1.682}}
			\multiput(67.339,31.534)(.0337004049,-.0361133603){247}{\line(0,-1){.0361133603}}
			\multiput(57.713,22.671)(.0364583333,.0336174242){264}{\line(1,0){.0364583333}}
			\multiput(75.331,8.39)(-.0419598109,.0336879433){423}{\line(-1,0){.0419598109}}
			\multiput(57.457,8.265)(.0428463357,.0336879433){423}{\line(1,0){.0428463357}}
			\put(57.723,22.73){\line(1,0){18.135}}
			\put(14.419,4){\makebox(0,0)[cc]{$A_{iii}$}}
			\put(40.432,4){\makebox(0,0)[cc]{$A_{iv}$}}
			\put(65.851,4){\makebox(0,0)[cc]{$A_v$}}
		\end{picture}
	
		\caption{Non-Isomorphic Simple Graphs on $5$ vertices and $7$ edges}
			\label{7}
		\end{figure}
	\end{center}
	The graph $A_v$ contains three 2-edge cuts. Therefore, adding a loop and a parallel edge or two parallel edges to $A_v$ still contains a 2-edge cut. Hence we discard $A_v$. If $G$ is obtained from $A_{iii}$ by adding a loop and a parallel edge, then $G$ contains a 2-edge cut. Suppose $G$ is obtained from $A_{iii}$ by adding two parallel edges such that no 2-edge cut is there in $G$. Then we get a contradiction by Lemma \ref{qlemma}. Thus $G$ can not be obtained from $A_{iii}$. If $G$ is obtained from $A_{iv}$ by adding a loop and a parallel edge or two parallel edges such that $G$ does not contain 2-edge, then $M(G) \cong {\tilde{M}({G_2})}$.  Thus $G$ can not be obtained from $A_{iv}$.

	{\bf Case-(iii)} If $G $ contains either a loop and two 2-cycles or three 2-cycles.

\noindent Then $G$ is obtained from a 2-connected simple graph on $5$ vertex and $6$ edges. By Harary (\cite{Har}, Pg 217), there are two non-isomorphic 2-connected simple graphs on $5$ vertices and $6$ edges as shown in Figure \ref{8}.
	 
	\begin{center}
		\begin{figure}[h]
		\unitlength 1mm 
		\linethickness{0.4pt}
		\ifx\plotpoint\undefined\newsavebox{\plotpoint}\fi 
		\begin{picture}(51.079,32.04)(0,0)
			\put(5.848,8.171){\circle*{1.682}}
			\put(23.717,8.138){\circle*{1.682}}
			\put(6.172,22.425){\circle*{1.682}}
			\put(24.041,22.393){\circle*{1.682}}
			\put(5.956,22.427){\line(0,-1){14.125}}
			\put(23.831,8.052){\line(0,1){14.625}}
			\put(6.058,8.171){\line(1,0){18.079}}
			\put(15.511,31.199){\circle*{1.682}}
			\multiput(15.562,31.399)(.0337004049,-.0361133603){247}{\line(0,-1){.0361133603}}
			\multiput(5.936,22.536)(.0364583333,.0336174242){264}{\line(1,0){.0364583333}}
			\put(32.045,8.064){\circle*{1.682}}
			\put(49.914,8.031){\circle*{1.682}}
			\put(32.369,22.318){\circle*{1.682}}
			\put(50.238,22.286){\circle*{1.682}}
			\put(32.153,22.32){\line(0,-1){14.125}}
			\put(50.028,7.945){\line(0,1){14.625}}
			\put(41.708,31.092){\circle*{1.682}}
			\multiput(41.759,31.292)(.0337004049,-.0361133603){247}{\line(0,-1){.0361133603}}
			\multiput(32.133,22.429)(.0364583333,.0336174242){264}{\line(1,0){.0364583333}}
			\put(5.946,22.595){\line(1,0){18.135}}
			\multiput(32.257,22.595)(.0409824052,-.0336889263){428}{\line(1,0){.0409824052}}
			\multiput(32.108,8.027)(.0432238999,.0337357268){423}{\line(1,0){.0432238999}}
			\put(14.27,4.459){\makebox(0,0)[cc]{$A_{vi}$}}
			\put(40.581,4.459){\makebox(0,0)[cc]{$A_{vii}$}}
		\end{picture}
			\caption{Non-Isomorphic 2-Connected Graphs on 5 Vertices and 6 Edges}
			\label{8}
			
		\end{figure}
	\end{center}
	If $G$ is obtained from $A_{vi}$ or $A_{vii}$ by adding a loop and two parallel edges, then $G$ contains a 2-edge cut, which is a contradiction by Lemma $\ref{glem}$. Suppose $G$ obtained from $A_{vi}$ by adding  three parallel edges such that $G$ does not contain a 2-edge cut.  Then $G \cong a_{i}$ or $G\cong a_{ii}$ where $a_{i}$, $a_{ii}$ are shown in Figure \ref{qdiscard}. By Lemma \ref{qlemma}, we discard both the graphs $a_{i}$ and $a_{ii}$. \\
\begin{figure}[h!]
\centering
\unitlength 1mm 
\linethickness{0.4pt}
\ifx\plotpoint\undefined\newsavebox{\plotpoint}\fi 
\begin{picture}(107.398,32.703)(0,0)
	\put(5.848,8.171){\circle*{1.682}}
	\put(23.717,8.138){\circle*{1.682}}
	\put(6.172,22.425){\circle*{1.682}}
	\put(24.041,22.393){\circle*{1.682}}
	\put(5.956,22.427){\line(0,-1){14.125}}
	\put(23.831,8.052){\line(0,1){14.625}}
	\put(6.058,8.171){\line(1,0){18.079}}
	\put(15.511,31.199){\circle*{1.682}}
	\multiput(15.562,31.399)(.0337004049,-.0361133603){247}{\line(0,-1){.0361133603}}
	\multiput(5.936,22.536)(.0364583333,.0336174242){264}{\line(1,0){.0364583333}}
	\put(32.045,8.064){\circle*{1.682}}
	\put(49.914,8.031){\circle*{1.682}}
	\put(32.369,22.318){\circle*{1.682}}
	\put(50.238,22.286){\circle*{1.682}}
	\put(32.153,22.32){\line(0,-1){14.125}}
	\put(50.028,7.945){\line(0,1){14.625}}
	\put(32.255,8.064){\line(1,0){18.079}}
	\put(41.708,31.092){\circle*{1.682}}
	\multiput(41.759,31.292)(.0337004049,-.0361133603){247}{\line(0,-1){.0361133603}}
	\multiput(32.133,22.429)(.0364583333,.0336174242){264}{\line(1,0){.0364583333}}
	\put(5.946,22.595){\line(1,0){18.135}}
	\put(32.108,22.446){\line(1,0){18.284}}
	\put(57.625,8.306){\circle*{1.682}}
	\put(75.494,8.273){\circle*{1.682}}
	\put(57.949,22.56){\circle*{1.682}}
	\put(75.818,22.528){\circle*{1.682}}
	\put(57.733,22.562){\line(0,-1){14.125}}
	\put(75.608,8.187){\line(0,1){14.625}}
	\put(67.288,31.334){\circle*{1.682}}
	\multiput(67.339,31.534)(.0337004049,-.0361133603){247}{\line(0,-1){.0361133603}}
	\multiput(57.713,22.671)(.0364583333,.0336174242){264}{\line(1,0){.0364583333}}
	\put(84.084,8.009){\circle*{1.682}}
	\put(101.953,7.976){\circle*{1.682}}
	\put(84.408,22.263){\circle*{1.682}}
	\put(102.277,22.231){\circle*{1.682}}
	\put(84.192,22.265){\line(0,-1){14.125}}
	\put(102.067,7.89){\line(0,1){14.625}}
	\put(93.747,31.037){\circle*{1.682}}
	\multiput(93.798,31.237)(.0337004049,-.0361133603){247}{\line(0,-1){.0361133603}}
	\multiput(84.172,22.374)(.0364583333,.0336174242){264}{\line(1,0){.0364583333}}
	\multiput(101.79,8.093)(-.0419598109,.0336879433){423}{\line(-1,0){.0419598109}}
	\multiput(83.916,7.968)(.0428463357,.0336879433){423}{\line(1,0){.0428463357}}
	\qbezier(5.946,8.324)(13.155,.595)(23.635,8.324)
	\qbezier(23.635,8.324)(30.027,16.5)(23.932,22.297)
	\qbezier(23.932,22.297)(24.973,29.878)(15.311,31.513)
	\qbezier(31.811,8.176)(41.473,1.115)(49.946,8.027)
	\qbezier(50.243,22.446)(56.932,14.79)(49.946,8.027)
	\qbezier(32.257,22)(31.067,30.473)(41.77,31.216)
	\qbezier(57.824,22.743)(59.162,32.703)(67.04,31.365)
	\qbezier(75.811,22.446)(81.979,14.493)(75.365,8.027)
	\qbezier(75.365,8.027)(67.635,.297)(57.527,8.027)
	\qbezier(84.135,22.149)(78.338,15.311)(83.837,7.878)
	\qbezier(102.27,22.446)(107.398,16.426)(102.121,7.73)
	\qbezier(93.946,31.067)(102.27,30.399)(102.27,22.297)
	\put(57.675,8.176){\line(1,0){17.838}}
	\qbezier(67.27,31.533)(75.995,32.164)(75.889,22.704)
	\put(13.454,1.472){\makebox(0,0)[cc]{$a_i$}}
	\put(40.783,1.472){\makebox(0,0)[cc]{$a_{ii}$}}
	\put(66.85,1.472){\makebox(0,0)[cc]{$a_{iii}$}}
	\put(92.287,1.472){\makebox(0,0)[cc]{$a_{iv}$}}
\end{picture}

\caption{Graphs Obtained from $A_{vi}$, $A_{vii}$ and From a cycle on 5 Vertices}
\label{qdiscard}
\end{figure}
	 Suppose $G$ is obtained from $A_{vii}$ by adding three parallel edges such that $G$ does not contain 2-edge cut then $G \cong H_{11}$ or $G\cong a_{iv}$, where $a_{iv}$ is shown in Figure \ref{qdiscard}. Suppose $G\cong a_{iv}$, since due to cycles of $G$, $N \backslash a$ does not contain circuit with 6 elements, thus $N \backslash a \ncong K_{3,3}$. Hence $G \ncong a_{iv}$.

	 	{\bf Case-(iv)} If $G$ contains either a loop and three pairs of parallel edges or four pairs of parallel edges.

\noindent Then $G$ is obtained from a 2-connected simple graph on $5$ vertices and $5$ edges by adding a loop and three parallel edges or four parallel edges. By Harary (\cite{Har}, Pg 216), a cycle on $5$ vertices is the only a 2-connected graph on  $5$ vertices and $5$ edges. Thus $G$ is obtained by adding a loop and three parallel edges such that $G$ does not contain a 2-edge cut, but there is no graph without a 2-edge cut. Suppose $G$ is obtained from a cycle on $5$ vertices by adding four parallel edges such that it does not contain a 2-edge cut. Then $a_{iii}$ is the only graph with this condition, where $a_{iii}$ is a graph shown in Figure \ref{qdiscard}. By Lemma \ref{qlemma}, $G\ncong a_{iii}$.

	{\bf Case (v)} If $G$ contains a loop and four pairs of parallel edges or five pairs of parallel edges.

\noindent  Then $G$ is obtained from a 2-connected simple graph on $5$ vertices and four edges. By Harary (\cite{Har}, Pg 216), there is no 2-connected simple graph on $5$ vertices and $4$ edges.

Therefore $G$ is isomorphic to $H_i$, for some $i \in \{7,8,9,10,11\}$.
\end{proof}

\section{Forbidden Minors for the class $\mathcal{M}_2$}
This section gives alternative proof for the forbidden minor characterization of the class $\mathcal{M}_2$ using the results obtained from the last two sections. Recall that $\mathcal{M}_k$ is the class of graphic matroid whose splitting with respect to any set of $k$-elements is cographic and $\mathcal{F}=\{F_7, F^*_7, M(K_5), M(K_{3,3})\}.$

We need the following lemmas, which follow easily from the definition of the splitting operation. 
\begin{lem}\label{splem}
	Let $M$ be a binary matroid such that it contains more than two 2-circuits. Then $M_T$ will contain at least one 2-circuit for any $T \subseteq E(M)$ with $|T|=2$. 
\end{lem}

\begin{lem}\label{looplem}
	If $M(G)$ is a graphic matroid and contains a loop, then $M(G)_T$ is graphic with $|T|=2$ or includes a loop. 
\end{lem}

\begin{lem}\label{paralem}
	Let M(G) be a graphic matroid such that it contains two adjacent 2-circuits. Then $M(G)_T$ is graphic for some $T\subseteq E(M)$ with $|T|=2$ or $M_T$ contains a 2-circuit. 
\end{lem}

Now, we give an alternate prove to Theorem \ref{gct}. Here we restate the result along with the graphs. 
\begin{thm}
	Let $M$ be a graphic matroid without containing a minor isomorphic to any member of $\mathcal{T}_2$. Then $M \in \mathcal{M}_2$ if and only if $M$ does not contain a minor isomorphic to any one of the circuit matroids $M(G_1)$, $M(G_2)$ and $M(G_3)$, where $G_1$, $G_2$ and $G_3$ are the graphs as shown in Figure \ref{M2}.
\end{thm}
\begin{figure}[h!]
	\centering
\unitlength 1mm 
\linethickness{0.4pt}
\ifx\plotpoint\undefined\newsavebox{\plotpoint}\fi 
\begin{picture}(89.119,31.693)(0,0)
	\put(3.569,8.056){\circle*{1.682}}
	\put(21.438,8.024){\circle*{1.682}}
	\put(3.893,22.31){\circle*{1.682}}
	\put(21.762,22.278){\circle*{1.682}}
	\put(3.751,22.316){\line(1,0){18}}
	\put(3.677,22.312){\line(0,-1){14.125}}
	\put(21.552,7.937){\line(0,1){14.625}}
	\qbezier(3.569,22.561)(12.924,29.078)(21.859,22.141)
	\put(3.779,8.056){\line(1,0){18.079}}
	\qbezier(3.359,8.056)(13.24,.067)(21.438,8.056)
	\multiput(3.674,22.351)(.0418961088,-.0337152118){424}{\line(1,0){.0418961088}}
	\multiput(3.359,8.056)(.0433835446,.0337152118){424}{\line(1,0){.0433835446}}
	\put(28.16,7.824){\circle*{1.682}}
	\put(46.029,7.791){\circle*{1.682}}
	\put(28.484,22.078){\circle*{1.682}}
	\put(46.353,22.046){\circle*{1.682}}
	\put(28.342,22.084){\line(1,0){18}}
	\put(28.268,22.08){\line(0,-1){14.125}}
	\put(46.143,7.705){\line(0,1){14.625}}
	\put(28.37,7.824){\line(1,0){18.079}}
	\qbezier(27.95,7.824)(37.83,-.165)(46.029,7.824)
	\put(37.823,30.852){\circle*{1.682}}
	\multiput(37.874,31.052)(.0337022512,-.0361095549){247}{\line(0,-1){.0361095549}}
	\multiput(28.248,22.189)(.0364583607,.0336174495){264}{\line(1,0){.0364583607}}
	\multiput(27.998,7.939)(.033670059,.0778620114){297}{\line(0,1){.0778620114}}
	\put(61.028,8.421){\circle*{1.682}}
	\put(78.897,8.388){\circle*{1.682}}
	\put(61.352,22.675){\circle*{1.682}}
	\put(79.221,22.643){\circle*{1.682}}
	\put(61.21,22.681){\line(1,0){18}}
	\put(79.011,8.302){\line(0,1){14.625}}
	\put(61.238,8.421){\line(1,0){18.079}}
	\multiput(61.133,22.716)(.0418962264,-.0337146226){424}{\line(1,0){.0418962264}}
	\multiput(60.818,8.421)(.0433820755,.0337146226){424}{\line(1,0){.0433820755}}
	\put(53.028,16.398){\circle*{1.682}}
	\put(88.278,16.523){\circle*{1.682}}
	\multiput(52.687,16.432)(8.93751,.03125){4}{\line(1,0){8.93751}}
	\multiput(79.062,8.432)(.038381772,.033713718){241}{\line(1,0){.038381772}}
	\multiput(88.312,16.557)(-.050403264,.033602176){186}{\line(-1,0){.050403264}}
	\multiput(52.937,16.682)(.047857179,.033571454){175}{\line(1,0){.047857179}}
	\multiput(52.812,16.307)(.035326113,-.033695677){230}{\line(1,0){.035326113}}
	\put(13.258,0.237){\makebox(0,0)[cc]{$G_1$}}
	\put(37.3,0.237){\makebox(0,0)[cc]{$G_2$}}
	\put(68.766,0.237){\makebox(0,0)[cc]{$G_3$}}
\end{picture}
	
	\caption{Forbidden Minors for the Class $\mathcal{M}_2$}
	\label{M2}
\end{figure}
\begin{proof}
	Let $M$ be a graphic matroid such that it contains a minor  isomorphic to   $M(G_i)$ for some $i \in \{1,2,3\}$. Then we can easily prove that $M_T$ is not cographic for some $T \subseteq E(M)$ and $|T|=2$.

 Conversely, assume that $M$ avoids  a minor  isomorphic to   $M(G_i)$ for every $i \in \{1,2,3\}$. Then we prove that $M \in \mathcal{M}_2.$  On contrary, assume that $M\notin \mathcal{M}_2$. Then, by Lemma \ref{cg},  $M$ contains a minor isomorphic to $F \in \mathcal{F}$.  Note that $\mathcal{M}_1$ is empty and $M$ does not contain a minor isomorphic to any element of $\mathcal{T}_2$. 
   Hence, by Lemma \ref{lm1}, $M$ has a minor $Q$ such that (i) $Q_T \cong F$, (ii) $Q_T / T'\cong F$ for $T' \subseteq T$. By Lemma \ref{lm2} and Lemma \ref{lm3},  $Q$ is without 2-cocircuit and a coloop. Therefore $Q$ is 2-minimal matroid with respect to $F \in \mathcal{F}$. 
  Thus, by Lemma \ref{ql},  there exist a binary matroid $N$ containing $a$ such that $N\backslash a \cong F$ and $Q=N/a$ or $Q$ is coextension of $N/a$ by one or two elements. Since $M$ is graphic, $Q$ is also graphic. Let $Q=M(G)$ for some connected graph $G$.  Since, by Definition \ref{spdef}, $Q_T$ contains a 2-cocircuit as $|T|=2$, $Q_T \ncong F$. Therefore $Q_T / T'\cong F.$ Hence $Q$ is a coextension of $N/a$ by one or two elements.

	\textbf{Case (i).} Suppose $F=F_7^*$.  Then $N\backslash a \cong F_7^*$ and , by Lemma $\ref{q1}$, $N/a \cong M(H_1)$ or $N/a \cong M(H_2),$ where $H_1$ and $H_2$ are the graphs shown in Figure \ref{3}. 
	  Therefore $G$ is a  coextension of $H_1$ or $H_2$ by one or two elements. 
	  Suppose $G$ is a coextension of $H_1$ by an element. 
	  By Lemma $\ref{glem}$ (4), $G$ can not contain a 2-edge cut. We take coextensions of $H_1$ such that it does not contain a 2-edge cut. There is only one non-isomorphic coextension $C_1$ of $H_1$ by an element as shown in Figure \ref{10}. Since  $C_1 \cong G_2$,  $G$ can not arise from $H_1$. 
	\begin{center}
		\begin{figure}[h]
		\unitlength 1mm 
		\linethickness{0.4pt}
		\ifx\plotpoint\undefined\newsavebox{\plotpoint}\fi 
		\begin{picture}(26.236,32.04)(0,0)
			\put(5.848,8.171){\circle*{1.682}}
			\put(23.717,8.138){\circle*{1.682}}
			\put(6.172,22.425){\circle*{1.682}}
			\put(24.041,22.393){\circle*{1.682}}
			\put(5.956,22.427){\line(0,-1){14.125}}
			\put(23.831,8.052){\line(0,1){14.625}}
			\put(15.511,31.199){\circle*{1.682}}
			\multiput(15.562,31.399)(.0337004049,-.0361133603){247}{\line(0,-1){.0361133603}}
			\multiput(5.936,22.536)(.0364583333,.0336174242){264}{\line(1,0){.0364583333}}
			\multiput(5.797,22.595)(.0420243308,-.0336889263){428}{\line(1,0){.0420243308}}
			\multiput(5.797,8.473)(.0436989436,.0336696779){415}{\line(1,0){.0436989436}}
			\qbezier(15.162,31.513)(26.236,30.919)(23.932,22.595)
			\put(15.162,4.459){\makebox(0,0)[cc]{$C_1$}}
			\put(5.797,8.027){\line(1,0){17.986}}
		\end{picture}
			\caption{Coextensions of $H_1$ by one element.}
			\label{10}
		\end{figure}
	\end{center}
	  
	Suppose $G$ arise from $H_2$.  Then $G$ is a coextension of $H_2$ by one or two elements. There are two non-isomorphic coextensions of $H_2$ by an element as shown in Figure \ref{11}, such that it does not contain a 2-edge cut. 
	\begin{center}
		\begin{figure}[h]
		\unitlength 1mm 
		\linethickness{0.4pt}
		\ifx\plotpoint\undefined\newsavebox{\plotpoint}\fi 
		\begin{picture}(51.079,32.257)(0,0)
			\put(5.848,8.171){\circle*{1.682}}
			\put(23.717,8.138){\circle*{1.682}}
			\put(6.172,22.425){\circle*{1.682}}
			\put(24.041,22.393){\circle*{1.682}}
			\put(5.956,22.427){\line(0,-1){14.125}}
			\put(23.831,8.052){\line(0,1){14.625}}
			\put(6.058,8.171){\line(1,0){18.079}}
			\put(15.511,31.199){\circle*{1.682}}
			\multiput(15.562,31.399)(.0337004049,-.0361133603){247}{\line(0,-1){.0361133603}}
			\multiput(5.936,22.536)(.0364583333,.0336174242){264}{\line(1,0){.0364583333}}
			\put(32.045,8.064){\circle*{1.682}}
			\put(49.914,8.031){\circle*{1.682}}
			\put(32.369,22.318){\circle*{1.682}}
			\put(50.238,22.286){\circle*{1.682}}
			\put(32.153,22.32){\line(0,-1){14.125}}
			\put(50.028,7.945){\line(0,1){14.625}}
			\put(41.708,31.092){\circle*{1.682}}
			\multiput(41.759,31.292)(.0337004049,-.0361133603){247}{\line(0,-1){.0361133603}}
			\multiput(32.133,22.429)(.0364583333,.0336174242){264}{\line(1,0){.0364583333}}
			\put(5.946,22.595){\line(1,0){18.135}}
			\put(32.257,22.297){\line(1,0){18.135}}
			\put(32.108,8.324){\line(1,0){17.838}}
			\multiput(49.797,8.473)(-.033726868,.09618403){238}{\line(0,1){.09618403}}
			\qbezier(5.797,22.595)(7.581,32.257)(15.311,31.216)
			\qbezier(24.081,22.743)(29.284,14.27)(23.784,8.176)
			\qbezier(5.649,8.176)(15.98,.669)(23.635,8.324)
			\qbezier(32.108,8.324)(41.919,.966)(49.946,8.176)
			\put(15.311,0.35){\makebox(0,0)[cc]{$C_2$}}
			\put(41.77,0.35){\makebox(0,0)[cc]{$C_3$}}
		\end{picture}
		
			\caption{Coextensions of $H_2$ by one element.}
			\label{11}

		\end{figure}
	\end{center} 
	 Here $C_{3} \cong G_2.$ Therefore $G \ncong C_3$. By Lemma $\ref{splem}$,  the splitting of $M(C_2)$ with respect to any two elements contains a 2-circuit.    Thus $(C_{2})_T \ncong F_7^*$ for any $T \in E(C_{2})$ with $|T|=2$ as $F_7^*$ is a 3-connected matroid.  Thus $G \ncong C_{2} $.  Also there is no coextension of $C_2$ by an element, which does not contain a 2-edge cut. Hence, we discard $C_2$.

	\textbf{{ Case (ii).}}  Suppose $F=F_7$. Then $N\backslash a \cong F_7$ and,  by Lemma \ref{q1}, $N/a \cong M(H_3),$ where $H_3$ is the graph as shown in Figure $\ref{3}$. Therefore $G$ is a coextension of $M(H_3)$ by one or two elements. Suppose $G$ is a coextension of $H_3$ by an element. There are five $C_4, C_5, C_6, C_7, C_8$ non-isomorphic coextensions of $M(H_3)$ by an element as shown in Figure $\ref{9}$. By Lemma \ref{splem}, the splitting of the matroids $M(C_4)$, $M(C_7)$ and $M(C_8)$  contains a 2-circuit after splitting with any two elements set. Also, there is no coextension of $C_4$, $C_7$ and $C_8$ by one element such that it does not contain a loop or 2-edge cut or a minor isomorphic to $M(G_1)$, $M(G_2)$ or $M(G_3)$. Hence we discard coextensions of $C_4$, $C_7$ and $C_8$. By Lemma \ref{looplem}, splitting of $M(C_5)$ with respect to any two elements is either graphic or contains a loop. Also coextension of $C_5$ by one element either contains a loop or a minor isomorphic to $M(G_2)$. Hence $G \ncong C_5$. 
		Here $M(C_6)\cong M(G_1) $. Hence $G\ncong C_6$. 
	
		\begin{figure}[h] 
\unitlength 1mm 
\linethickness{0.4pt}
\ifx\plotpoint\undefined\newsavebox{\plotpoint}\fi 
\begin{picture}(145.788,31.92)(0,0)
	\put(8.615,9.528){\circle*{1.682}}
	\put(26.484,9.496){\circle*{1.682}}
	\put(8.939,23.782){\circle*{1.682}}
	\put(26.808,23.75){\circle*{1.682}}
	\put(8.797,23.788){\line(1,0){18}}
	\put(8.723,23.784){\line(0,-1){14.125}}
	\put(26.598,9.409){\line(0,1){14.625}}
	\put(8.825,9.528){\line(1,0){18.079}}
	\put(36.894,9.473){\circle*{1.682}}
	\put(54.763,9.44){\circle*{1.682}}
	\put(37.218,23.727){\circle*{1.682}}
	\put(55.087,23.695){\circle*{1.682}}
	\put(37.002,23.729){\line(0,-1){14.125}}
	\put(54.877,9.354){\line(0,1){14.625}}
	\put(37.104,9.473){\line(1,0){18.079}}
	\qbezier(36.684,9.473)(46.564,1.484)(54.763,9.473)
	\qbezier(26.613,23.922)(20.425,31.524)(29.795,29.933)
	\qbezier(29.795,29.933)(34.744,25.779)(26.966,24.099)
	\put(64.061,9.25){\circle*{1.682}}
	\put(81.93,9.217){\circle*{1.682}}
	\put(64.385,23.504){\circle*{1.682}}
	\put(82.254,23.472){\circle*{1.682}}
	\put(64.243,23.51){\line(1,0){18}}
	\put(64.169,23.506){\line(0,-1){14.125}}
	\put(82.044,9.131){\line(0,1){14.625}}
	\put(64.271,9.25){\line(1,0){18.079}}
	\qbezier(63.851,9.25)(73.731,1.261)(81.93,9.25)
	\qbezier(64.362,23.699)(73.554,31.92)(82.039,23.523)
	\put(92.23,9.04){\circle*{1.682}}
	\put(110.099,9.007){\circle*{1.682}}
	\put(92.554,23.294){\circle*{1.682}}
	\put(110.423,23.262){\circle*{1.682}}
	\put(92.412,23.3){\line(1,0){18}}
	\put(92.338,23.296){\line(0,-1){14.125}}
	\put(110.213,8.921){\line(0,1){14.625}}
	\put(92.44,9.04){\line(1,0){18.079}}
	\qbezier(92.02,9.04)(101.9,1.051)(110.099,9.04)
	\qbezier(92.531,23.313)(84.929,18.009)(92.177,9.171)
	\qbezier(92.531,23.489)(101.723,31.71)(110.208,23.313)
	\put(120.4,8.83){\circle*{1.682}}
	\put(138.269,8.797){\circle*{1.682}}
	\put(120.724,23.084){\circle*{1.682}}
	\put(138.593,23.052){\circle*{1.682}}
	\put(120.582,23.09){\line(1,0){18}}
	\put(120.508,23.086){\line(0,-1){14.125}}
	\put(138.383,8.711){\line(0,1){14.625}}
	\put(120.61,8.83){\line(1,0){18.079}}
	\qbezier(120.19,8.83)(130.07,.841)(138.269,8.83)
	\qbezier(120.701,23.279)(129.893,31.5)(138.378,23.103)
	\multiput(36.999,23.755)(.0422440191,-.0336961722){418}{\line(1,0){.0422440191}}
	\multiput(36.999,9.67)(.043881068,.0336771845){412}{\line(1,0){.043881068}}
	\multiput(63.907,23.755)(.0421438679,-.0337146226){424}{\line(1,0){.0421438679}}
	\multiput(63.907,9.25)(.0442583732,.0336961722){418}{\line(1,0){.0442583732}}
	\multiput(120.246,8.829)(.0431367925,.0337146226){424}{\line(1,0){.0431367925}}
	\qbezier(26.698,23.965)(17.238,31.113)(8.619,23.965)
	\qbezier(8.409,9.67)(1.892,19.866)(8.829,23.755)
	\qbezier(26.067,9.67)(34.476,15.556)(26.908,23.965)
	\qbezier(54.868,9.46)(60.439,13.875)(59.703,6.517)
	\qbezier(59.703,6.937)(54.973,-.105)(54.868,9.25)
	\qbezier(110.156,23.335)(117.724,17.869)(110.156,9.04)
	\qbezier(138.325,8.829)(145.788,15.031)(138.536,23.335)
	\put(18.499,1){\makebox(0,0)[cc]{$C_4$}}
	\put(45.618,1){\makebox(0,0)[cc]{$C_5$}}
	\put(72.947,1){\makebox(0,0)[cc]{$C_6$}}
	\put(100.696,1){\makebox(0,0)[cc]{$C_7$}}
	\put(129.706,1){\makebox(0,0)[cc]{$C_8$}}
	\put(36.999,23.965){\line(1,0){18.289}}
\end{picture}

		\caption{Coextensions of $H_3$ by one Element}
		\label{9}
	\end{figure}

	\textbf{{ Case (iii).}} Suppose  $F=M(K_5).$ Then $N\backslash a \cong M(K_5)$ and, by Lemma \ref{q3}, $N / a \cong M(H_4) $, $N / a \cong M(H_5)$ or  $N / a \cong M(H_6)$, where $H_4$, $H_5$ and $H_6$ are the graphs shown in Figure \ref{4}.  Suppose $N / a \cong M(H_4)$. Then $M(H_4)_T \ncong M(K_5)$ for any $T$ with $|T|=2$. Thus we take coextension of $H_4$. Here all coextensions of $H_4$ by an element belongs to $\mathcal{T}_2$, hence we discard $H_4$.
 Suppose $N / a \cong M(H_5)$. Therefore $G$ is a coextension of $H_5$ by one or two elements. There are seven coextensions of $H_5$ by an element shown in Figure \ref{16}, without containing a minor isomorphic to $M(G_1)$, $M(G_2)$ and $M(G_3)$.   Since $M(C_i)$ contains a loop and at least two 2-circuits, the splitting with respect to any two elements set $T$ will contain a loop or a 2-circuit for $i=9,10, \cdots, 15$. Therefore we take coextensions of $C_i$ for $i=9,10, \cdots 15$ such that it does not contain a minor isomorphic to  $M(G_1)$ or $M(G_2)$ or $M(G_3)$. Also,  it does not contain a loop or more than two 2-circuits.  It is observed that there is no coextension with these conditions. Hence we discard $H_5$. Suppose $N / a \cong M(H_6)$. Then, by using similar argument as above, we  discard $H_6$. 
	\begin{figure}[h!]
		\centering
	\unitlength 1mm 
	\linethickness{0.4pt}
	\ifx\plotpoint\undefined\newsavebox{\plotpoint}\fi 
	\begin{picture}(106.94,69.841)(0,0)
		\put(5.848,8.171){\circle*{1.682}}
		\put(23.717,8.138){\circle*{1.682}}
		\put(6.172,22.425){\circle*{1.682}}
		\put(24.041,22.393){\circle*{1.682}}
		\put(5.956,22.427){\line(0,-1){14.125}}
		\put(23.831,8.052){\line(0,1){14.625}}
		\put(6.058,8.171){\line(1,0){18.079}}
		\put(15.511,31.199){\circle*{1.682}}
		\multiput(15.562,31.399)(.0337004049,-.0361133603){247}{\line(0,-1){.0361133603}}
		\multiput(5.936,22.536)(.0364583333,.0336174242){264}{\line(1,0){.0364583333}}
		\put(32.045,8.064){\circle*{1.682}}
		\put(49.914,8.031){\circle*{1.682}}
		\put(32.369,22.318){\circle*{1.682}}
		\put(50.238,22.286){\circle*{1.682}}
		\put(32.153,22.32){\line(0,-1){14.125}}
		\put(50.028,7.945){\line(0,1){14.625}}
		\put(32.255,8.064){\line(1,0){18.079}}
		\put(41.708,31.092){\circle*{1.682}}
		\multiput(41.759,31.292)(.0337004049,-.0361133603){247}{\line(0,-1){.0361133603}}
		\multiput(32.133,22.429)(.0364583333,.0336174242){264}{\line(1,0){.0364583333}}
		\put(5.946,22.595){\line(1,0){18.135}}
		\put(57.625,8.306){\circle*{1.682}}
		\put(75.494,8.273){\circle*{1.682}}
		\put(57.949,22.56){\circle*{1.682}}
		\put(75.818,22.528){\circle*{1.682}}
		\put(57.733,22.562){\line(0,-1){14.125}}
		\put(75.608,8.187){\line(0,1){14.625}}
		\put(67.288,31.334){\circle*{1.682}}
		\multiput(67.339,31.534)(.0337004049,-.0361133603){247}{\line(0,-1){.0361133603}}
		\multiput(57.713,22.671)(.0364583333,.0336174242){264}{\line(1,0){.0364583333}}
		\qbezier(31.811,8.176)(41.473,1.115)(49.946,8.027)
		\qbezier(75.365,8.027)(67.635,.297)(57.527,8.027)
		\put(57.675,8.176){\line(1,0){17.838}}
		\put(6.937,44.452){\circle*{1.682}}
		\put(24.806,44.419){\circle*{1.682}}
		\put(7.261,58.706){\circle*{1.682}}
		\put(25.13,58.674){\circle*{1.682}}
		\put(7.045,58.708){\line(0,-1){14.125}}
		\put(24.92,44.333){\line(0,1){14.625}}
		\put(7.147,44.452){\line(1,0){18.079}}
		\put(16.6,67.48){\circle*{1.682}}
		\multiput(16.651,67.68)(.0337004049,-.0361133603){247}{\line(0,-1){.0361133603}}
		\multiput(7.025,58.817)(.0364583333,.0336174242){264}{\line(1,0){.0364583333}}
		\put(33.134,44.345){\circle*{1.682}}
		\put(51.003,44.312){\circle*{1.682}}
		\put(33.458,58.599){\circle*{1.682}}
		\put(51.327,58.567){\circle*{1.682}}
		\put(33.242,58.601){\line(0,-1){14.125}}
		\put(33.344,44.345){\line(1,0){18.079}}
		\put(42.797,67.373){\circle*{1.682}}
		\multiput(42.848,67.573)(.0337004049,-.0361133603){247}{\line(0,-1){.0361133603}}
		\multiput(33.222,58.71)(.0364583333,.0336174242){264}{\line(1,0){.0364583333}}
		\put(58.714,44.587){\circle*{1.682}}
		\put(76.583,44.554){\circle*{1.682}}
		\put(59.038,58.841){\circle*{1.682}}
		\put(76.907,58.809){\circle*{1.682}}
		\put(58.822,58.843){\line(0,-1){14.125}}
		\put(76.697,44.468){\line(0,1){14.625}}
		\put(68.377,67.615){\circle*{1.682}}
		\multiput(68.428,67.815)(.0337004049,-.0361133603){247}{\line(0,-1){.0361133603}}
		\multiput(58.802,58.952)(.0364583333,.0336174242){264}{\line(1,0){.0364583333}}
		\put(87.906,45.972){\circle*{1.682}}
		\put(105.775,45.939){\circle*{1.682}}
		\put(88.23,60.226){\circle*{1.682}}
		\put(106.099,60.194){\circle*{1.682}}
		\put(88.014,60.228){\line(0,-1){14.125}}
		\put(105.889,45.853){\line(0,1){14.625}}
		\put(97.569,69){\circle*{1.682}}
		\multiput(97.62,69.2)(.0337004049,-.0361133603){247}{\line(0,-1){.0361133603}}
		\multiput(87.994,60.337)(.0364583333,.0336174242){264}{\line(1,0){.0364583333}}
		\multiput(105.612,46.056)(-.0419598109,.0336879433){423}{\line(-1,0){.0419598109}}
		\multiput(87.738,45.931)(.0428463357,.0336879433){423}{\line(1,0){.0428463357}}
		\qbezier(7.035,44.605)(14.244,36.876)(24.724,44.605)
		\qbezier(32.9,44.457)(42.562,37.396)(51.035,44.308)
		\qbezier(76.9,58.727)(83.068,50.774)(76.454,44.308)
		\put(58.764,44.457){\line(1,0){17.838}}
		\multiput(6.937,44.567)(.0336939489,.0791075323){287}{\line(0,1){.0791075323}}
		\multiput(6.727,44.567)(.0431349575,.0337146794){424}{\line(1,0){.0431349575}}
		\multiput(7.148,59.072)(.041151741,-.0337146794){424}{\line(1,0){.041151741}}
		\multiput(33.005,44.357)(.0426391534,.0337146794){424}{\line(1,0){.0426391534}}
		\multiput(42.675,67.481)(.0336353508,-.0924972146){250}{\line(0,-1){.0924972146}}
		\multiput(51.084,44.357)(-.0426391534,.0337146794){424}{\line(-1,0){.0426391534}}
		\put(58.862,58.862){\line(1,0){17.869}}
		\multiput(68.112,67.691)(.033739164,-.094296637){243}{\line(0,-1){.094296637}}
		\multiput(58.652,44.567)(.0431349575,.0337146794){424}{\line(1,0){.0431349575}}
		\put(87.662,46.038){\line(1,0){18.079}}
		\multiput(105.741,46.038)(-.0336353508,.0933380983){250}{\line(0,1){.0933380983}}
		\multiput(5.676,22.704)(.0422453568,-.0336957013){418}{\line(1,0){.0422453568}}
		\multiput(23.335,8.619)(-.033706311,.096683893){237}{\line(0,1){.096683893}}
		\multiput(15.346,31.533)(-.0336939489,-.0813049637){287}{\line(0,-1){.0813049637}}
		\multiput(5.676,8.199)(.0421433493,.0337146794){424}{\line(1,0){.0421433493}}
		\multiput(31.743,8.199)(.0337214481,.0782050604){293}{\line(0,1){.0782050604}}
		\multiput(32.164,8.409)(.0426391534,.0337146794){424}{\line(1,0){.0426391534}}
		\multiput(32.374,22.283)(.0418400904,-.0336761704){412}{\line(1,0){.0418400904}}
		\multiput(57.601,8.409)(.0336939489,.0805724866){287}{\line(0,1){.0805724866}}
		\multiput(57.811,22.914)(.040889657,-.0336738351){437}{\line(1,0){.040889657}}
		\multiput(57.39,8.199)(.0425330744,.0337331279){430}{\line(1,0){.0425330744}}
		\qbezier(6.727,44.567)(.631,53.817)(7.148,58.862)
		\qbezier(7.148,58.862)(.315,59.282)(4.835,64.748)
		\qbezier(4.835,64.748)(9.04,66.325)(7.358,59.072)
		\qbezier(33.215,44.357)(26.383,52.66)(33.425,58.862)
		\qbezier(33.425,58.862)(26.488,58.547)(30.482,64.538)
		\qbezier(30.482,64.538)(35.738,66.114)(33.425,58.862)
		\qbezier(33.005,44.567)(50.348,48.351)(51.294,58.862)
		\qbezier(58.231,44.987)(74.418,47.51)(77.151,58.441)
		\qbezier(76.31,44.777)(72.421,37.945)(79.043,40.783)
		\qbezier(79.043,40.783)(81.986,45.198)(76.52,44.567)
		\qbezier(87.662,45.618)(91.026,40.468)(84.299,42.465)
		\qbezier(84.299,42.465)(81.881,45.513)(87.452,46.459)
		\qbezier(87.452,46.459)(95.966,37.314)(105.741,46.249)
		\qbezier(88.083,60.544)(94.074,46.669)(105.531,46.249)
		\qbezier(57.601,8.199)(73.262,11.247)(75.89,23.124)
		\qbezier(66.85,31.743)(61.7,34.581)(67.06,36.158)
		\qbezier(67.06,36.158)(71.896,35.527)(67.481,31.533)
		\qbezier(50.243,22.704)(48.456,13.244)(31.954,7.988)
		\qbezier(32.164,22.704)(26.173,21.443)(29.01,26.067)
		\qbezier(29.01,26.067)(31.008,28.59)(32.584,22.704)
		\qbezier(23.965,8.199)(28.59,10.196)(27.329,5.045)
		\qbezier(27.329,5.045)(22.389,.841)(23.755,8.409)
		\put(15.346,38){\makebox(0,0)[cc]{$C_9$}}
		\put(42.044,38){\makebox(0,0)[cc]{$C_{10}$}}
		\put(67.691,38){\makebox(0,0)[cc]{$C_{11}$}}
		\put(96.491,38){\makebox(0,0)[cc]{$C_{12}$}}
		\put(14.926,1.472){\makebox(0,0)[cc]{$C_{13}$}}
		\put(41.624,1.472){\makebox(0,0)[cc]{$C_{14}$}}
		\put(67.06,1.472){\makebox(0,0)[cc]{$C_{15}$}}
	\end{picture}
		
		\caption{Coextensions of $H_5$ by one element}
		\label{16}
	\end{figure}
	
	\textbf{{Case (iv).}} Suppose $F=M(K_{3,3})$. Then, by Lemma \ref{q4}, $N / a \cong M(H_i)$ for $i=7,8,9,10,11$, where $H_i$ is the graph as shown in Figure \ref{5}, for $i=7,8,9,10,11$.

	 Suppose $N/a \cong M(H_7)$, then $M(H_7)_T\ncong M(K_{3,3})$ for any $T$ with $|T|=2$, thus we take coextension of $H_7$, but all coextension of $H_7$ belongs to $\mathcal{T}_2$. Hence we discard the graph $H_7$.
	 Let $N/a \cong M(H_8)$. Then $G$ is a coextension of $H_8$ by one or two  elements. There are three coextensions of $H_8$ by an element as shown in Figure \ref{14}, without containing a 2-edge cut. Here $M(C_{16})$ and $M(C_{17})$ contains a minor isomorphic to $M(G_1)$ and $M(C_{18})$ contains a minor isomorphic to $M(G_3)$. Thus we discard these coextensions. 
	\begin{figure}[h!] 
		\centering
	\unitlength 1mm 
	\linethickness{0.4pt}
	\ifx\plotpoint\undefined\newsavebox{\plotpoint}\fi 
	\begin{picture}(122.572,22.596)(0,0)
		\put(53.67,7.37){\circle*{1.682}}
		\put(71.539,7.337){\circle*{1.682}}
		\put(53.994,21.624){\circle*{1.682}}
		\put(71.863,21.592){\circle*{1.682}}
		\put(53.852,21.63){\line(1,0){18}}
		\put(53.88,7.37){\line(1,0){18.079}}
		\multiput(53.775,21.665)(.0418962264,-.0337146226){424}{\line(1,0){.0418962264}}
		\multiput(53.46,7.37)(.0433820755,.0337146226){424}{\line(1,0){.0433820755}}
		\put(45.67,15.347){\circle*{1.682}}
		\put(80.92,15.472){\circle*{1.682}}
		\multiput(71.704,7.381)(.038381743,.033713693){241}{\line(1,0){.038381743}}
		\multiput(80.954,15.506)(-.050403226,.033602151){186}{\line(-1,0){.050403226}}
		\multiput(45.579,15.631)(.047857143,.033571429){175}{\line(1,0){.047857143}}
		\multiput(45.454,15.256)(.035326087,-.033695652){230}{\line(1,0){.035326087}}
		\put(12.2,7.238){\circle*{1.682}}
		\put(30.069,7.205){\circle*{1.682}}
		\put(12.524,21.492){\circle*{1.682}}
		\put(30.393,21.46){\circle*{1.682}}
		\put(12.382,21.498){\line(1,0){18}}
		\put(30.183,7.119){\line(0,1){14.625}}
		\put(12.41,7.238){\line(1,0){18.079}}
		\put(4.2,15.215){\circle*{1.682}}
		\put(39.45,15.34){\circle*{1.682}}
		\multiput(30.234,7.249)(.038381743,.033713693){241}{\line(1,0){.038381743}}
		\multiput(39.484,15.374)(-.050403226,.033602151){186}{\line(-1,0){.050403226}}
		\multiput(4.109,15.499)(.047857143,.033571429){175}{\line(1,0){.047857143}}
		\multiput(3.984,15.124)(.035326087,-.033695652){230}{\line(1,0){.035326087}}
		\put(94.481,7.501){\circle*{1.682}}
		\put(112.35,7.468){\circle*{1.682}}
		\put(94.805,21.755){\circle*{1.682}}
		\put(112.674,21.723){\circle*{1.682}}
		\put(94.663,21.761){\line(1,0){18}}
		\put(94.691,7.501){\line(1,0){18.079}}
		\multiput(94.586,21.796)(.0418962264,-.0337146226){424}{\line(1,0){.0418962264}}
		\put(86.481,15.478){\circle*{1.682}}
		\put(121.731,15.603){\circle*{1.682}}
		\multiput(112.515,7.512)(.038381743,.033713693){241}{\line(1,0){.038381743}}
		\multiput(121.765,15.637)(-.050403226,.033602151){186}{\line(-1,0){.050403226}}
		\multiput(86.265,15.387)(.035326087,-.033695652){230}{\line(1,0){.035326087}}
		\multiput(3.994,15.346)(.113987259,-.033635585){225}{\line(1,0){.113987259}}
		\multiput(11.983,7.568)(.115311867,.033706546){237}{\line(1,0){.115311867}}
		\put(12.403,21.863){\line(0,-1){14.716}}
		\multiput(45.408,15.556)(.108139097,-.033739398){243}{\line(1,0){.108139097}}
		\multiput(80.725,15.556)(-.114424853,-.033706546){237}{\line(-1,0){.114424853}}
		\multiput(85.981,15.346)(.142771365,.033725519){187}{\line(1,0){.142771365}}
		\multiput(85.981,15.767)(.107273984,-.033739398){243}{\line(1,0){.107273984}}
		\multiput(94.39,7.568)(.1093156499,.0336355846){250}{\line(1,0){.1093156499}}
		\put(94.6,21.863){\line(0,-1){14.505}}
		\put(21.443,2.523){\makebox(0,0)[cc]{$C_{16}$}}
		\put(62.226,2.523){\makebox(0,0)[cc]{$C_{17}$}}
		\put(102.799,2.523){\makebox(0,0)[cc]{$C_{18}$}}
	\end{picture}
	
		\caption{Coextensions of $H_8$ by one element}
\label{14}
\end{figure}

Let $N/a \cong M(H_9)$.  Thus we take the coextension of $H_9$ by one or two elements. There is no coextension of $H_9$ by one element as well as by two elements, such that it does not contain a loop, a 2-cocircuit or a minor isomorphic to $M(G_1)$, $M(G_2)$ or $M(G_3)$. Hence we discard $H_9$.

	   If  $N/a \cong M(H_{10})$ then, $G$ is a coextension of $H_{10}$ by one  or two elements such that it does not contain 2-cocircuit or minor  $M(G_1)$, $M(G_2)$ or $M(G_3)$.  Such non-isomorphic coextensions of $H_{10}$ by one element are shown in Figure \ref{15}.
	\begin{figure}[h!]
		\centering
	\unitlength 1mm 
	\linethickness{0.4pt}
	\ifx\plotpoint\undefined\newsavebox{\plotpoint}\fi 
	\begin{picture}(94.389,23.545)(0,0)
		\put(12.2,7.238){\circle*{1.682}}
		\put(30.069,7.205){\circle*{1.682}}
		\put(12.524,21.492){\circle*{1.682}}
		\put(30.393,21.46){\circle*{1.682}}
		\put(12.382,21.498){\line(1,0){18}}
		\put(12.41,7.238){\line(1,0){18.079}}
		\put(4.2,15.215){\circle*{1.682}}
		\put(39.45,15.34){\circle*{1.682}}
		\multiput(30.234,7.249)(.038381743,.033713693){241}{\line(1,0){.038381743}}
		\multiput(39.484,15.374)(-.050403226,.033602151){186}{\line(-1,0){.050403226}}
		\multiput(4.109,15.499)(.047857143,.033571429){175}{\line(1,0){.047857143}}
		\multiput(3.984,15.124)(.035326087,-.033695652){230}{\line(1,0){.035326087}}
		\put(61.687,7.501){\circle*{1.682}}
		\put(79.556,7.468){\circle*{1.682}}
		\put(62.011,21.755){\circle*{1.682}}
		\put(79.88,21.723){\circle*{1.682}}
		\put(61.869,21.761){\line(1,0){18}}
		\put(61.897,7.501){\line(1,0){18.079}}
		\put(53.687,15.478){\circle*{1.682}}
		\put(88.937,15.603){\circle*{1.682}}
		\multiput(79.721,7.512)(.038381743,.033713693){241}{\line(1,0){.038381743}}
		\multiput(53.471,15.387)(.035326087,-.033695652){230}{\line(1,0){.035326087}}
		\multiput(12.193,7.148)(.0431349575,.0337146794){424}{\line(1,0){.0431349575}}
		\put(12.193,21.653){\line(5,-4){17.869}}
		\put(3.784,15.136){\line(1,0){35.527}}
		\multiput(53.396,15.767)(.046457667,.033681808){181}{\line(1,0){.046457667}}
		\multiput(61.595,21.863)(.145018725,-.033725285){187}{\line(1,0){.145018725}}
		\multiput(88.713,15.556)(-5.015271,.030032){7}{\line(-1,0){5.015271}}
		\multiput(61.385,7.568)(.0436307616,.0337146794){424}{\line(1,0){.0436307616}}
		\put(79.884,22.073){\line(0,-1){14.715}}
		\qbezier(39.311,15.346)(42.254,22.494)(45.198,15.346)
		\qbezier(45.198,15.346)(43.411,8.199)(39.522,15.346)
		\qbezier(88.713,15.977)(92.392,23.545)(93.969,15.977)
		\qbezier(93.969,15.977)(94.389,10.301)(88.923,15.556)
		\put(20.812,2.943){\makebox(0,0)[cc]{$C_{19}$}}
		\put(70.424,2.943){\makebox(0,0)[cc]{$C_{20}$}}
	\end{picture}
		
		\caption{Coextensions of $H_{10}$ by one element}
		\label{15}
	\end{figure}
	Here $C_{19}$ and $C_{20}$ contains a loop.  Since $M(C_{19}) \in \mathcal{T}_2$, $G \ncong C_{19}$. By Lemma \ref{looplem},  the splitting of $M(C_{20})$ contains a loop or its graphic. Hence we take coextensions by $H_{10}$ by two elements such that it does not contain a loop, a 2-cocircuit and a minor isomorphic to $M(G_1)$, $M(G_2)$ or $M(G_3)$. It is observed that there is no such coextension. Thus we discard $H_{10}$.

	 If $N/a \cong M(H_{11})$, then,  we take coextensions of $H_{11}$ by one element such that it does not contain 2-cocircuit and minor $M(G_1)$, $M(G_2)$ or $M(G_3)$.  One can easily check that there is no such coextension. Hence we discard $H_{11}$.
	
	Therefore $M_T$ is cographic for any $T$ with $|T|=2$. Hence proved. \end{proof}

\section{Forbidden-Minor Characterization of the class $\mathcal{M}_3$ }
In this section, we prove Theorem \ref{mt}. We now first prove 'if part' of Theorem \ref{mt} in the following proposition.
\begin{prop}\label{prop}
	Let $M$ be a graphic matroid. If $M$ contains a minor isomorphic to any of the matroids  $M(\tilde{G_1})$, $M(\tilde{G_3})$ and $M(G_i)$ for $i=4,5,6,7$, then $M_T$ is not cographic for some $T\subseteq E(M)$ with $|T|=3$, where $G_i$ is the graph as shown in Figure \ref{M3}, for $i=1,3,4,5,6,7.$ \end{prop}
\begin{proof}
	Note that, by Theorem \ref{gct}, $M(G_1)$, $M(G_2)$ and $M(G_3)$ are the forbidden minors for the class $\mathcal{M}_2$. For $k=1,3$, let $\tilde{M}({G_k})$ be a matroid obtained by adding an element $z_k$ to $M(G_k)$. 	Then $\tilde{M}({G_k})\backslash z_k = M(G_k)$. By Theorem \ref{gct}, there is $T_k$ with $|T_k|=2$ such that $M(G_k)_{T_k}$ is not a cographic matroid.  Let $T=T_k \cup z_k$.  Then, by Lemma \ref{bp1}(ii), $\tilde{M}({G_k})_T\backslash z_k \cong (\tilde{M}({G_k})\backslash z_k)_{T_k} \cong M(G_k)_{T_k}\cong M(G_k)_{T_k}$ is not cographic. $\tilde{M}({G_k})_T$ is not cographic. Also $G_2$ is minor of $G_5$.

	Suppose $M=M(G_4)$ and $T=\{x,y,z\}$, where $x,y,z$ are as shown in the graph $G_4$. Let $A_4$ be the standard matrix representation of $G_4$ as follows:\\
	
	$ A_4 =\left[ \begin{array}{ccccccc}
		x & y & z & & & & \\
		1 & 0 & 0 &0&1&0&1 \\
		0 & 1 & 0 &0&1&1&0 \\
		0 & 0 & 0 &1&1&1&1 \end{array} \right].$ Then $(A_4)_T=\left[ \begin{array}{ccccccc}
		x & y & z & & & & \\
		1 & 0 & 0 &0&1&0&1 \\
		0 & 1 & 0 &0&1&1&0 \\
		0 & 0 & 0 &1&1&1&1 \\
		1 & 1 & 1 &0&0&0&0\end{array} \right]$ \\
	Let $(A_4)_T^{'}$ be the matrix obtained by applying the row transformations $R_4 \implies R_4+(R_1+R_2)$ and $R_{34}$. Then we have\\
	\begin{center}
	$(A_4)_T^{'} =\left[ \begin{array}{ccccccc}
		1 & 0 & 0 &0&1&0&1 \\
		0 & 1 & 0 &0&1&1&0 \\
		0 & 0 & 1 &0&0&1&1\\
		0 & 0 & 0 &1&1&1&1 \end{array} \right]$ \end{center}
	
	Since the vector matroid of $(A_4)_T^{'}$ is isomorphic to $F_7^*$, 
	 $M(G_4)_T \cong M[(A_4)_T] \cong M[(A_4)_T^{'}] \cong F_7^* $. By Theorem \ref{cg}, $M(G_4)_T$ is not cographic. Similarly, we can prove that $M(G_5)_{T_5} \cong F_7^*$, $M(G_6)_{T_6} \cong M(K_{3,3})$ and $M(G_{7})_{T_7} \cong M(K_{3,3})$ for $T_i=\{x,y,z\}$ , where $x,y,z$ are as shown in the graph $G_i$ in Figure \ref{M3}, for $i=5, 6, 7$.

	Suppose $M$ contains a minor $N$ isomorphic to any of circuit matroids $M(G_k)$, $k=4,5,6,7$. Then $N_T$ is not cographic for some $T \subseteq E(N)$ with $|T|=3$. Since $N$ is a minor of $M$, there exist subsets $T_1$ and $T_2$ of $E(M)$ such that $M\backslash T_1 /T_2 \cong N$. Then, by repeated application of Lemma \ref{bp1}, we have $M_T\backslash T_1 /T_2 \cong (M\backslash T_1/T_2)_T \cong N_T$. As $N_T $ is not cographic, $M_T$ is also not cographic.
\end{proof}
Now we prove the main theorem. For convenience, we restate the main theorem here along with graphs.

\begin{thm}
Let $M$ be a graphic matroid without containing a minor isomorphic to any member of the set $\mathcal{T}_3.$ Then $M \in \mathcal{M}_3$ if and only if $M$ does not contain a minor isomorphic to any of the matroid $\tilde{M}(G_{1})$, $\tilde{M}(G_{3})$, $M(G_4)$,  $M(G_5)$, $M(G_6)$ and $M(G_7)$,  where $G_i$ is the graph as shown in Figure \ref{M3}, for $i=1,3,4,5,6,7$.
\end{thm}
\begin{figure}[h!]
	\centering
\unitlength 1mm 
\linethickness{0.4pt}
\ifx\plotpoint\undefined\newsavebox{\plotpoint}\fi 
\begin{picture}(107.374,59.935)(0,0)
	\put(3.569,8.056){\circle*{1.682}}
	\put(21.438,8.024){\circle*{1.682}}
	\put(3.893,22.31){\circle*{1.682}}
	\put(21.762,22.278){\circle*{1.682}}
	\put(3.751,22.316){\line(1,0){18}}
	\put(3.677,22.312){\line(0,-1){14.125}}
	\put(21.552,7.937){\line(0,1){14.625}}
	\put(3.779,8.056){\line(1,0){18.079}}
	\multiput(3.674,22.351)(.0418962264,-.0337146226){424}{\line(1,0){.0418962264}}
	\multiput(3.359,8.056)(.043384434,.0337146226){424}{\line(1,0){.043384434}}
	\put(33.11,8.001){\circle*{1.682}}
	\put(50.979,7.968){\circle*{1.682}}
	\put(33.434,22.255){\circle*{1.682}}
	\put(51.303,22.223){\circle*{1.682}}
	\put(33.292,22.261){\line(1,0){18}}
	\put(33.218,22.257){\line(0,-1){14.125}}
	\put(51.093,7.882){\line(0,1){14.625}}
	\put(33.32,8.001){\line(1,0){18.079}}
	\qbezier(32.9,8.001)(42.78,.012)(50.979,8.001)
	\qbezier(33.411,22.274)(25.809,16.97)(33.057,8.132)
	\qbezier(33.411,22.45)(42.603,30.671)(51.088,22.274)
	\qbezier(21.567,22.45)(15.379,30.052)(24.749,28.461)
	\qbezier(24.749,28.461)(29.698,24.307)(21.92,22.627)
	\put(58.293,8.166){\circle*{1.682}}
	\put(76.162,8.133){\circle*{1.682}}
	\put(58.617,22.42){\circle*{1.682}}
	\put(76.486,22.388){\circle*{1.682}}
	\put(58.401,22.422){\line(0,-1){14.125}}
	\put(76.276,8.047){\line(0,1){14.625}}
	\put(58.503,8.166){\line(1,0){18.079}}
	\put(67.956,31.194){\circle*{1.682}}
	\multiput(68.007,31.394)(.0337004049,-.0361133603){247}{\line(0,-1){.0361133603}}
	\multiput(58.381,22.531)(.0364583333,.0336174242){264}{\line(1,0){.0364583333}}
	\multiput(58.131,8.281)(.0336700337,.0778619529){297}{\line(0,1){.0778619529}}
	\multiput(67.999,31.375)(.033613445,-.097163866){238}{\line(0,-1){.097163866}}
	\multiput(75.999,8.25)(-.0419598109,.0336879433){423}{\line(-1,0){.0419598109}}
	\multiput(58.125,8.125)(.0428463357,.0336879433){423}{\line(1,0){.0428463357}}
	\put(83.34,8.864){\circle*{1.682}}
	\put(101.209,8.831){\circle*{1.682}}
	\put(83.664,23.118){\circle*{1.682}}
	\put(101.533,23.086){\circle*{1.682}}
	\put(101.323,8.745){\line(0,1){14.625}}
	\put(83.55,8.864){\line(1,0){18.079}}
	\qbezier(83.13,8.864)(93.01,.875)(101.209,8.864)
	\put(93.003,31.892){\circle*{1.682}}
	\multiput(93.054,32.092)(.0337004049,-.0361133603){247}{\line(0,-1){.0361133603}}
	\multiput(83.428,23.229)(.0364583333,.0336174242){264}{\line(1,0){.0364583333}}
	\multiput(83.178,8.979)(.0336700337,.0778619529){297}{\line(0,1){.0778619529}}
	\multiput(83.499,23.125)(.0426610979,-.0337112172){419}{\line(1,0){.0426610979}}
	\qbezier(101.374,23)(107.374,17.375)(101.374,8.75)
	\qbezier(83.499,23.25)(96.374,21.875)(101.249,9)
	\put(12.193,1){\makebox(0,0)[cc]{$G_4$}}
	\put(42.255,1){\makebox(0,0)[cc]{$G_5$}}
	\put(67.271,1){\makebox(0,0)[cc]{$G_6$}}
	\put(93.759,1){\makebox(0,0)[cc]{$G_7$}}
	\put(5.256,38.913){\circle*{1.682}}
	\put(23.125,38.881){\circle*{1.682}}
	\put(5.58,53.167){\circle*{1.682}}
	\put(23.449,53.135){\circle*{1.682}}
	\put(5.438,53.173){\line(1,0){18}}
	\put(5.364,53.169){\line(0,-1){14.125}}
	\put(23.239,38.794){\line(0,1){14.625}}
	\qbezier(5.256,53.418)(14.611,59.935)(23.546,52.998)
	\put(5.466,38.913){\line(1,0){18.079}}
	\qbezier(5.046,38.913)(14.927,30.924)(23.125,38.913)
	\multiput(5.361,53.208)(.0418962264,-.0337146226){424}{\line(1,0){.0418962264}}
	\multiput(5.046,38.913)(.043384434,.0337146226){424}{\line(1,0){.043384434}}
	\put(39.591,39.278){\circle*{1.682}}
	\put(57.46,39.245){\circle*{1.682}}
	\put(39.915,53.532){\circle*{1.682}}
	\put(57.784,53.5){\circle*{1.682}}
	\put(39.773,53.538){\line(1,0){18}}
	\put(57.574,39.159){\line(0,1){14.625}}
	\put(39.801,39.278){\line(1,0){18.079}}
	\multiput(39.696,53.573)(.0418962264,-.0337146226){424}{\line(1,0){.0418962264}}
	\multiput(39.381,39.278)(.0433820755,.0337146226){424}{\line(1,0){.0433820755}}
	\put(31.591,47.255){\circle*{1.682}}
	\put(66.841,47.38){\circle*{1.682}}
	\multiput(31.25,47.289)(8.9375,.03125){4}{\line(1,0){8.9375}}
	\multiput(57.625,39.289)(.038381743,.033713693){241}{\line(1,0){.038381743}}
	\put(66.875,47.414){\line(-3,2){9.375}}
	\multiput(31.5,47.539)(.047857143,.033571429){175}{\line(1,0){.047857143}}
	\multiput(31.375,47.164)(.035326087,-.033695652){230}{\line(1,0){.035326087}}
	\put(14.945,32.094){\makebox(0,0)[cc]{$G_1$}}
	\put(47.329,32.094){\makebox(0,0)[cc]{$G_3$}}
	\put(42.254,27.749){\makebox(0,0)[cc]{$x$}}
	\put(41.624,5.886){\makebox(0,0)[cc]{$z$}}
	\put(34.897,15.556){\makebox(0,0)[cc]{$y$}}
	\put(61.595,27.959){\makebox(0,0)[cc]{$x$}}
	\put(74.418,27.749){\makebox(0,0)[cc]{$y$}}
	\put(67.06,9.46){\makebox(0,0)[cc]{$z$}}
	\put(92.077,6.727){\makebox(0,0)[cc]{$z$}}
	\put(94.599,20.602){\makebox(0,0)[cc]{$y$}}
	\put(106.162,16.607){\makebox(0,0)[cc]{$x$}}
	\put(24.175,26.488){\makebox(0,0)[cc]{$z$}}
	\put(12.403,24.596){\makebox(0,0)[cc]{$x$}}
	\put(11.983,9.46){\makebox(0,0)[cc]{$y$}}
\end{picture}
	
	\caption{Forbidden Minor for the Class $\mathcal{M}_3$}
	\label{M3}
\end{figure}
\begin{proof}
	If $M$ has a minor isomorphic to one of the circuit matroids $\tilde{M}({G_1})$, and $\tilde{M}({G_3})$ and $M(G_k)$ for $k=4,5,6,7$, then $M_T$ is not cographic for some $T \in E(N)$ with $|T|=3$ by Proposition \ref{prop}.

Conversely, assume that $M$ does not contain a minor isomorphic to any of the circuit matroid $\tilde{M}({G_1}),$  $\tilde{M}({G_3})$ and $M(G_k)$ for $k=4,5,6,7$. Then we will prove that $M_T$ is cographic for any $T \subseteq E(M)$ with $|T|=3$. On contrary, assume that $M_T$ is not cographic for some $T \subseteq E(M)$ with $|T|=3$.  Then, by Theorem \ref{cg} $M_T$, has a minor isomorphic to $F \in \mathcal{F}$.
Then, by Lemma \ref{lm1}, $M$ has a minor $Q$ containing T such that one of the following holds.
	\begin{enumerate}[(i).]
	\item $Q$ isomorphic to one of the circuit matroids $\tilde{M}({G_1}),$ $\tilde{M}({G_2})$ and $\tilde{M}({G_3}).$  
	\item $Q_T \cong F$ or $Q_T/T^{'} \cong F$ for some $T^{'} \subseteq T$.
	\end{enumerate}
	
	Note that the graph $G_2$ contains a minor isomorphic to the graph  $G_5$. 
	If $Q$ is isomorphic to  $\tilde{M}({G_1}),$ $\tilde{M}({G_2})$ and $\tilde{M}({G_3}),$ then  it is clear that $M$ contains a minor isomorphic to $\tilde{M}({G_1})$, $\tilde{M}({G_3})$ or $M(G_5)$, a contradiction. Thus (i) does not hold.
	Suppose (ii) holds. By Lemma \ref{lm2} and Lemma \ref{lm3}, assume that $Q$ do not contain a 2-cocircuit and a coloop. Thus $Q$ is 3-minimal with respect to $F$. Then, by Lemma \ref{ql}, there exist a binary matroid $N$ such that $N\backslash a \cong F$ and $N/a=Q$ or $Q$ is a coextension of $N/a$ by one, two or three elements. If $M$ is graphic, then $Q$ and $N/a$ are also graphic.  Let $Q\cong M(G)$ for some connected graph $G.$  We assume that $G$ does not contain a loop or 2-edge cut.

	\textbf{Case (i).}  Suppose $F=F_7^*$.  Then $N \backslash a = F_7^*$ and, by Lemma \ref{q1}, $N/a \cong M(H_1) $ or $M(H_2),$ where $H_1$ and $H_2$ are the graphs as shown in Figure \ref{3}. Here $H_1 \cong G_4$ and $H_2 \cong G_5$. Thus $Q=N/a \cong M(G_4)$ or $ Q=N/a \cong M(G_5)$. Hence $M$ contains a minor isomorphic to $M(G_4)$ or $M(G_5)$, a contradiction.

		\textbf{Case (ii).} Suppose $F=F_7$. Then $N\backslash a= F_7$ and, by Lemma \ref{q1}, $N/a \cong M(H_{3})$, where $H_{3}$ is the graph as shown in Figure \ref{3}. The graph $H_{3}$ contain three pairs parallel edges and a loop. Therefore $M(H_{3})_T$ contains at least one 2-circuit or a loop for any $T$ with $|T|=3$. Therefore $M(H_{3})_T \ncong F_7$ and $Q \ncong M(H_{3})$. Hence $Q$ is a coextension of $M(H_{3})$ by one, two or three elements. There are five non-isomorphic coextension of $M(H_{3})$ as shown in Figure \ref{9}, without containing 2-edge cut.     	The graphs $C_{4},$ $C_7$ and $C_8$ contains a minor isomorphic to the graph $G_5$ and $C_{5}$ contains a minor isomorphic to graph $G_4$.   Therefore we discard them.

		 Let $Q=M(C_{6})$ or $Q$ is a coextension of $M(C_{6})$ by one or two elements. Coextensions of $M(C_{6})$ by one or two elements contains 2-cocircuit or $G_5$ as a minor, a contradiction. Hence $Q=M(C_{6})$. Since $C_{6}$ has two 2-circuits and eight edges whereas $F_7$ has seven edges, there is $T \subseteq E(C_{6})$ with $|T|=3$ such that $M(C_{6})_T/x \cong F_7$ for some $x \in T $. If a 2-circuit does not intersect $T$ or contained in $T$, then it remains a 2-circuit in $M(C_{6})_T$. Thus $M(C_{6})_T/x \cong F_7$ contains a 2-circuit or a loop, a contradiction.  Hence $T$ contains exactly one edge from each 2-circuits of $C_6$. We assume that $x \in T$ does not belong to any  2-circuit of $C_6$. Then there is triangle in $G$ containing two edges from $T$ one of which is $x$. Then that 3-circuit of $M(C_{6})$ is get preserved in $M(C_6)_T$. Therefore $Q_T/x \cong F_7$ contains a 2-circuit, a contradiction.

	\textbf{Case (iii).} Suppose $F= M(K_5)$. Then $N\backslash a \cong M(K_5)$ and, by Lemma \ref{q3}, $N/a \cong M(H_4)$, $N/a \cong M(H_5)$ or $N/a \cong M(H_6)$, where $H_4$, $H_5$  and $H_6$ are the graphs as shown in Figure \ref{4}.  The matroid $M(H_4)$ does not give $M(K_5)$ after splitting with respect to a three elements set. Hence we take coextension of $M(H_4)$, but coextension of $M(H_4)$ is a trivial minor for the class $\mathcal{M}_2$. Hence we discard $M(H_4)$. The matroid $M(H_5)$ contains a minor isomorphic to $M(G_4)$. Hence $M$ can not be arise from  $M(H_5)$.  The matroid $M(H_6)$ contains a minor isomorphic to $M(G_4)$, and hence $G$ not arise from $H_5$ and $H_6$.

	\textbf{Case (iv).} Suppose $F= M(K_{3,3})$. Then $N\backslash a \cong M(K_{3,3})$ and, by Lemma \ref{q4}, $N/a \cong M(H_i)$,   where $H_i$ is the graph shown in Figure \ref{5}, for $i=7,8,9,10,11$.  The matroid $M(H_7)$ belongs to $\mathcal{T}_2$. Hence we discard the graph $H_7$. 
Note that $H_8 \cong G_6$ and $H_{11}=G_7$. The matroid $M(H_{i})$, for $i=9,10$ contains a minor isomorphic to $M(G_4)$. Therefore $M$ can not be arise from these matroid.

Thus from all the above cases, $M_T$ does not contain $F$ as a minor for any $F \in \mathcal{F}$. Therefore $M_T$ is cographic and hence proved. \end{proof}

\bibliographystyle{amsplain}

\end{document}